\documentclass{article}

\usepackage[utf8]{inputenc} 
\usepackage[T1]{fontenc}    
\usepackage[hidelinks]{hyperref}       
\usepackage{url}            
\usepackage{booktabs}       
\usepackage{amsfonts}       
\usepackage{nicefrac}       
\usepackage{microtype}      
\usepackage{xcolor}         

\usepackage{xr}
\usepackage{amssymb}
\usepackage{amsmath}
\usepackage{amstext}
\usepackage{amsthm}
\usepackage[linesnumbered,ruled,vlined]{algorithm2e}
\usepackage{algorithmic}
\usepackage{enumitem}
\usepackage{fullwidth}
\usepackage{float}
\usepackage{colonequals}
\usepackage{graphicx}
\usepackage{caption}
\usepackage{subcaption}
\usepackage{authblk}
\usepackage{comment}

\newtheorem{theorem}{Theorem}[section]
\newtheorem{lemma}[theorem]{Lemma}

\newtheorem{proposition}[theorem]{Proposition}
\newtheorem{corollary}[theorem]{Corollary}

\newtheorem{remark}[theorem]{Remark}

\newcommand{\BB}{\mathbb{B}}
\newcommand{\RR}{\mathbb{R}}

\newcommand{\NN}{\mathbb{N}}

\newcommand{\PP}{\mathbb{P}}
\newcommand{\EE}{\mathbb{E}}
\newcommand{\GG}{\mathbb{G}}

\renewcommand{\SS}{\mathbb{S}}

\newcommand{\sA}{\mathcal{A}}

\newcommand{\sC}{\mathcal{C}}

\newcommand{\sF}{\mathcal{F}}

\newcommand{\sH}{\mathcal{H}}

\newcommand{\sP}{\mathcal{P}}

\usepackage[capitalize]{cleveref}
\usepackage[notcite,notref,final]{showkeys}

\newcommand*{\rd}{\mathrm{d}}
\newcommand*{\dd}{\, \rd}

\title{Distributional Convergence of the Sliced Wasserstein Process}

\author[1]{Jiaqi Xi}
\author[1,2]{Jonathan Niles-Weed \thanks{JX and JNW are supported in part by NSF Grant DMS-2015291.}}
\affil[1]{Courant Institute of Mathematical Sciences, New York University, NY 10012}
\affil[2]{Center for Data Science, New York University, NY 10011}

\date{}

\begin{document}

\maketitle

\begin{abstract}
Motivated by the statistical and computational challenges of computing Wasserstein distances in high-dimensional contexts, machine learning researchers have defined modified Wasserstein distances based on computing distances between one-dimensional projections of the measures.
Different choices of how to aggregate these projected distances (averaging, random sampling, maximizing) give rise to different distances, requiring different statistical analyses.
We define the \emph{Sliced Wasserstein Process}, a stochastic process defined by the empirical Wasserstein distance between projections of empirical probability measures to all one-dimensional subspaces, and prove a uniform distributional limit theorem for this process.
As a result, we obtain a unified framework in which to prove distributional limit results for all Wasserstein distances based on one-dimensional projections.
We illustrate these results on a number of examples where no distributional limits were previously known.
\end{abstract}

\section{Introduction}
The Wasserstein distances have become useful tools in machine learning and data science, with applications in transfer learning~\cite{CouFlaTui17,RedHabSeb17}, generative modeling~\cite{GenPeyCut18,bousquet2017optimal}, statistics~\cite{GhoSen19,DebGhoSen21}, and various scientific domains~\cite{schiebinger2019optimal,dai2018autoencoder}.
Despite the popularity of these distances, they suffer from serious drawbacks in high dimensions.
From a statistical standpoint, estimating the Wasserstein distances from data suffers from the \emph{curse of dimensionality}, with convergence rates degrading sharply as the dimension increases~\cite{Dud69,weedBach,SinPoc18,nr-stm}.
From a computational standpoint, despite recent algorithmic advances~\cite{Cut13,AltWeeRig17}, the best algorithms for approximately computing general Wasserstein distances between distributions supported on $n$ points in $\RR^d$ for $d \geq 2$ have running times scaling quadratically in $n$, which is prohibitive on very large data sets.
These deficiencies have motivated the development of modifications of the Wasserstein distances which reduce the high-dimensional case to a series of one-dimensional problems.

Explicitly, given two compactly supported probability distributions $P$ and $Q$ in $\RR^d$, we write $P_u$ and $Q_u$ for the projections of $P$ and $Q$ onto the one-dimensional subspace spanned by $u$, for any $u \in \SS^{d-1}$.
Explicitly, if $X \sim P$, we let $P_u$ denote the law of $u^\top X$.
The measures $P_u$ and $Q_u$ are probability distributions on $\RR$ obtained by collapsing $P$ and $Q$ to the one-dimensional ``slice" in the direction of $u$.
Crucially, no matter how large $d$ is, the Wasserstein distance $W_p^p(P_u, Q_u)$ between the one dimensional measures is always easy to work with: it can be estimated from data at the rate that is independent of the dimension, and if $P_u$ and $Q_u$ are supported on $n$ points, then $W_p^p(P_u, Q_u)$ can be computed in \emph{nearly linear} time by a simple sorting procedure.

This observation has given rise to a number of different proposals for defining a distance between $P$ and $Q$ by aggregating the one-dimensional distances, the most prominent of which is the sliced Wasserstein distance~\cite{rpdb-sliced,BonRabPey15}:
\begin{equation}\label{sw}
	SW_p^p(P, Q) := \int W_p^p(P_u, Q_u) \dd \sigma(u)\,,
\end{equation}
where $\sigma$ denotes the uniform measure on $\SS^{d-1}$. Other options include:\\
$\bullet$ Discrete Sliced Wasserstein distance: $\widehat{SW}_p^p(P, Q) := \frac 1 L \sum_{i=1}^L W_p^p(P_{u_i}, Q_{u_i})$,  $\{u_i\} \subseteq \SS^{d-1}$~\cite{BonRabPey15}.\\
$\bullet$ Max-Sliced Wasserstein distance: $MSW_p^p(P, Q) := \max_{u \in \SS^{d-1}} W_p^p(P_u, Q_u)$~\cite{nr-stm,deshpande2019max,paty2019subspace}.\\
$\bullet$ Distributional Sliced Wasserstein distance: $DSW_p^p(P, Q) := \sup_{\tau \in \sP_C} \int W_p^p(P_u, Q_u) \dd \tau(u)$, where $\sP_C$ is a subset of probability measures on $\SS^{d-1}$~\cite{nguyen2020distributional}. \\
\\
Though the details of these techniques differ, they can be put on a common footing: if we view the function $W: u \mapsto W_p^p(P_u, Q_u)$ as a bounded function on $\SS^{d-1}$, then each of these sliced distances takes the form $F(W)$ for some function $F: \ell^\infty(\SS^{d-1}) \to \RR$.

Since these distances are all based on one-dimensional projections, it is natural to conjecture that they enjoy improved statistical performance.
This conjecture has been verified in certain special cases~\cite{nr-stm,nadjahi2020statistical,mbw-trim}
However, the analysis of these distances has largely been conducted separately, with different arguments tailored to each distance.
This raises the following fundamental question: is there a \emph{unified} approach to the analysis of these distances, which provides statistical guarantees for all of them simultaneously?

In this work, we develop such a unified approach.
In addition to generalizing prior works, our techniques allow us to prove new distributional convergence results for the sliced Wasserstein distance and its many variants.
These results make it possible to construct asymptotically valid confidence intervals for variants of the sliced Wasserstein distances and to guarantee the validity of the bootstrap.
Prior to our work, such results were only known for the standard sliced Wasserstein distance~\eqref{sw}~\cite{mbw-trim} or for the sliced and max-sliced Wasserstein distances between discrete distributions~\cite{okano2022inference}.\footnote{Concurrently and independently of our work, two other groups obtained similar results:~\cite{goldfeld2022statistical,XuHua22} proved distributional limits for the sliced and max-sliced Wasserstein, but not for other variants. \cite{goldfeld2022statistical} also obtain much more general distributional limit results for other quantities related to the Wasserstein distances.}

Obtaining distributional limits for empirical Wasserstein distances is an active area of research.
In the one-dimensional case, fundamental contributions were made by~\cite{del1999central,del2005asymptotics}, and further progress has been made in the case where one or both of the measures are discrete~\cite{sommerfeld2018inference,tameling2019empirical,del2022central}.
Multi-dimensional limits were recently obtained by~\cite{bl-clt,del2021central}, but these are not centered at the population-level quantities, making them of limited utility for inference. However, when the distributions are very smooth, there exist estimators with distributional limits with good centering~\cite{manole2021plugin}.
In this work, we draw on techniques recently proposed in~\cite{hundrieser2022unifying} to obtain central limit theorems by exploiting duality.

We consider compactly supported probability measures $P$ and $Q$ in $\RR^d$ with connected support, and the Wasserstein distances $W_p^p$ for $p > 1$.
To analyze the empirical behavior of the sliced Wasserstein distance and its variants, we define a stochastic process
\begin{equation}\label{swp}
	\GG_n(u) := \sqrt{n}(W_p^p(P_{n u}, Q_{n u}) - W_p^p(P_u, Q_u))\, \quad \quad u \in \SS^{d-1}\,.
\end{equation}
where $P_n$ and $Q_n$ consist of i.i.d.\ samples.
We may view $\GG_n$ as a random element of $\ell^{\infty}(\SS^{d-1})$, which records the deviation of the Wasserstein distance from its population counterpart along every direction simultaneously.
We call $\GG_n$ the \emph{Sliced Wasserstein Process}.

Our main result shows that
\begin{equation}\label{swp}
 \GG_n \rightsquigarrow \GG \quad \text{in } \ell^\infty(\SS^{d-1})\,,
\end{equation}
where $\GG$ is a tight Gaussian process on $\SS^{d-1}$.
That is, the collection of random variables $\sqrt{n}(W_p^p(P_{n u}, Q_{n u}) - W_p^p(P_u, Q_u))$ indexed by elements of $\SS^{d-1}$ enjoys a \emph{uniform} central limit theorem.
As is well known in the statistics literature~\cite{vw-ep}, uniform central limit theorems of this type give rise to distributional limits for any sufficiently regular functional on $\ell^{\infty}(\SS^{d-1})$ via the functional delta method---in particular, we directly obtain distributional limit theorems for the sliced Wasserstein distance and its many variants as a special case.
Our results likewise give techniques for proving the consistency of the bootstrap for any of the mentioned functionals as a consequence of general results for uniform central limit theorems.

\section{Main Result} \label{result}
Throughout, $P$ and $Q$ denote two probability distributions in $\RR^d$ with compact support, contained in a closed ball $\overline{B(0,R)}$ around the origin.
We fix a $p > 1$, and consider the Wasserstein distance of order $p$:
\begin{equation}\label{primal}
	W_p^p(P, Q) = \inf_{\pi \in \Pi(P, Q)} \int \|x - y\|^p \dd \pi(x, y)\,,
\end{equation}
where the infimum is taken over all couplings between $P$ and $Q$.
It is well known (see \cite{vil-ot}) that this problem possesses a dual formulation:
\begin{equation}\label{dual}
	W_p^p(P, Q) = \sup_{f: \overline{B(0,R)} \to \RR} \int f \dd P + \int f^c \dd Q\,,
\end{equation}
where $f^c$ denotes the $c$-transform: $f^c(y) = \inf_{x \in \overline{B(0,R)}} \|x - y\|^p - f(x)$.
It can be shown~(e.g., \cite[Lemma 1 \& 5]{manole2021sharp}) that the supremum in this dual formulation is achieved, and that without loss of generality we may assume that $f$ satisfies $f(0) = 0$ and $\|f\|_{\mathrm{Lip}} \leq p (2R)^{p-1}$.
We denote the class of such functions by $\sC$, and call any maximizer a \emph{Kantorovich potential}.
In order to obtain Gaussian limits,we adopt the following crucial assumption:
\begin{enumerate}[label=(\textbf{CC})]
	\item \label{CC} For all $u \in \SS^{d-1}$, the support of $P_u$ or $Q_u$ is an interval.
\end{enumerate}
Assumption \ref{CC} guarantees that the supremum in~\eqref{dual} is achieved at a \emph{unique} Kantorovich potential in $\sC$~\cite[Proposition 7.18]{San15}. In fact, it is known that the conclusion can still hold under weaker conditions, but we do not pursue this refinement here~\cite{sha-kanto,yang2021optimal}.
In the absence of this uniquness, Gaussian limits fail to hold for the optimal transport problem, even for discrete measures~\cite{sommerfeld2018inference}.

We now state our main result.
\begin{theorem} \label{thm:1}
    Suppose that $P$ and $Q$ are two probability distributions in $\RR^d$ whose supports are contained in the closed $d$-ball $\overline{B(0,R)}$ for some $R > 0$. Assume that $P$ and $Q$ satisfy \ref{CC}.
    Let $P_n$ and $Q_m$ denote empirical measures consisting of $n$ and $m$ i.i.d.\ samples from $P$ and $Q$, respectively.
    If $n/(n+m) \to \lambda \in (0, 1)$ as $n, m \to \infty$, then
    \begin{equation} \label{eq:5}
	    \sqrt{\frac{nm}{n+m}}\left(W_p^p(P_{n\cdot},Q_{m\cdot}) - W_p^p(P_{\cdot},Q_{\cdot})\right) \rightsquigarrow \GG \quad \text{in } \ell^{\infty}(\SS^{d-1})\,,
    \end{equation}
	where $\GG$ is a tight zero-mean Gaussian process on $\SS^{d-1}$ with covariance function
	\begin{equation} \label{eq:6}
 	\begin{split} \EE\GG(u)\GG(v) =  & (1-\lambda) \int f_u(u^\top x) f_v(v^\top x) \dd P(x) + \lambda \int f_u^c(u^\top y) f_v^c(v^\top y) \dd Q(y) \\
	& - (1-\lambda) \int f_u(u^\top x) \dd P(x) \int f_v(v^\top x) \dd P(x) \\
	& - \lambda  \int f^c_u(u^\top y) \dd Q(x) \int f_v^c(v^\top Y) \dd Q(x) , \end{split}
    \end{equation} where $f_u, f_v \in \sC$ are the unique Kantorovich potentials for $(P_u, Q_u)$ and $(P_v, Q_v)$, respectively.
\end{theorem}
\Cref{thm:1} formally includes the one-sample case as well, by taking $\lambda = 0, 1$.

As alluded to above, \cref{thm:1} gives rise to a wealth of statistical theorems as easy corollaries.
To describe these implications, we return to the abstract setting described above: denote by $W: \SS^{d-1} \to \RR$ the function $W(u) = W_p^p(P_u, Q_u)$, and consider any functional $F: \ell^{\infty}(\SS^{d-1}) \to \RR$. Then the distances we consider take the form $F(W)$. By different choices of $F$, we obtain the sliced Wasserstein distance, the max-sliced Wasserstein distance, and the other variants described above. \Cref{thm:1} will allow us to compare $F(W)$ to its empirical counterpart $F(W_{nm})$, where $W_{nm}(u) = W_p^p(P_{nu}, Q_{m u})$.

We recall the definition of directional Hadamard differentiability \cite{rom-delta}: we say that $F$ is directionally Hadamard differentiable at $\Phi$ if for all sequences $h_n \searrow 0$ and $\Psi_n \to \Psi \in \ell^\infty(\SS^{d-1})$, the limit
\begin{equation*}
	\lim_{n \to \infty} \frac{F(\Phi + h_n \Psi_n) - F(\Phi)}{h_n} =: F'_{\Phi}(\Psi)
\end{equation*}
exists.
We verify in \cref{applications} the directional Hadamard differentiability of several examples.
Under this assumption, we have the following.
\begin{corollary}
	Assume $F$ is directionally Hadamard differentiable.
	Under the assumptions of \cref{thm:1},
	\begin{equation*}
		\sqrt{\frac{nm}{n+m}}(F(W_{nm}) - F(W)) \rightsquigarrow F'_W(\GG)\,.
	\end{equation*}
\end{corollary}
\begin{proof}
	See \cite{rom-delta}.
\end{proof}

We also obtain a consistency result for the bootstrap, which we state for simplicity in the $n = m$ case.

\begin{corollary}
	Assume that $F$ is directionally Hadamard differentiable, and adopt the assumptions of \cref{thm:1}.
	For $k \ll n$ denote by $P^*$ and $Q^*$ bootstrap empirical measures consisting of $k$ i.i.d.\ draws from $P_n$ and $Q_m$, respectively, and set $W^* = W_p^p(P^*_u, Q^*_u)$. If $k \to \infty$ and $k/n \to 0$, then
	\begin{multline*}
		\sup_{h \in BL(1)} \EE[h(\sqrt{k}(F(W^*))- F(W_n))|X_1, \dots, X_n, Y_1, \dots, Y_n] \\ - \EE[h(\sqrt{n}(F(W_{n}) - F(W)))] \overset{p}{\to} 0\,,
	\end{multline*}
	where $BL(1)$ is the set of functions with bounded Lipschitz norm 1.
\end{corollary}
\begin{proof}
	See \cite{dumbgen1993nondifferentiable}.
\end{proof}

Finally, when the functional $F: \ell^\infty(\SS^{d-1}) \to \RR$ has a \emph{linear} Hadamard derivative, the resulting statistic will again be asymptotically Gaussian.
The following uniform convergence result shows that we can consistently estimate the covariance function of $\GG$ from data, which can be used to obtain asymptotic confidence intervals in this setting.
\begin{theorem} \label{thm:2}
    Under the same assumptions as Theorem \ref{thm:1}, there exists an estimator $\{\hat{\Sigma}_{u,v}\}_{u,v \in \SS^{d-1}}$ for the covariance functions $\{\Sigma_{u,v}\}_{u,v \in \SS^{d-1}}$ of the limiting process $\GG$ in the sense that \begin{equation} \label{eq:21}
        \EE_{P,Q}\mathop{\sup}\limits_{u,v \in \SS^{d-1}}\|\hat{\Sigma}_{u,v} - \Sigma_{u,v}\|_{\infty} \to 0 \quad \text{ as } n \to 0.
    \end{equation}
\end{theorem}

\section{Applications} \label{applications}
In this section, we focus on the two most popular of the variants we discussed---the sliced and max-sliced Wasserstein distances---and show how our main results obtained in the previous section can be used to obtain accurate asymptotic inference for these quantities.
For notational simplicity, we focus on the case where $n = m$, and rescale the resuling Gaussian process by a factor of $\sqrt 2$.
\subsection{Sliced Wasserstein Distance} \label{appsw}
Asymptotic and finite-sample inference for the sliced Wasserstein distance (SW) has already been thoroughly studied by~\cite{mbw-trim}.
We show that we can recover some of their results from our techniques.
Their main focus was on a robustification of the SW distance, the trimmed SW distance, defined as
\begin{equation*}
    SW_{p,\delta}(P,Q) := \left(\int_{\SS^{d-1}}\int_{\delta}^{1-\delta}|F_{\theta}^{-1}(t) - G_{\theta}^{-1}(t)|^{p}\,dt\,d\sigma^d(\theta)\right)^{\frac{1}{p}},
\end{equation*} where $\sigma^d$ denotes the uniform probability measure on $\SS^{d-1}$ and $F_{\theta}^{-1},\,G_{\theta}^{-1}$ are the (pseudo-)inverse of the CDFs of $P_{\theta},\,Q_{\theta}$ respectively. When $\delta = 0$, $SW_{p,\delta}$ reduces to the original sliced Wasserstein distance $SW_p$.
This robustification is necessary when $P$ and $Q$ are no longer assumed to have compact support.
\cite{mbw-trim} derive Gaussian limits and boostrap consistency for this functional.

To see how these results (under our stricter assumptions) can also be obtained from \cref{thm:1}, we denote  by $F: \ell^\infty(\SS^{d-1}) \to \RR$ the integration functional:
\begin{equation*}
	F(\Phi) = \int \Phi(u) \dd \sigma(u)\,.
\end{equation*}
The dominated convergence theorem immediately implies that $F$ is Hadamard differentiable, with derivative
\begin{equation*}
	F'_\Phi(\Psi) =  \int \Psi(u) \dd \sigma(u)\,.
\end{equation*}

We obtain the following.
\begin{theorem} \label{thm:3}
     Suppose that two compactly supported probability distributions $P$ and $Q$ in $\RR^d$ satisfy \ref{CC}.
     Then
      \begin{equation} \label{eq:29}
         \sqrt{n}(SW_p^p(P_n,Q_n) - SW_p^p(P,Q)) \stackrel{\text{d}}{\to} S:= \int_{\SS^{d-1}}\GG(\theta)\,d\sigma^d(\theta).
     \end{equation}
\end{theorem}
The random variable $S$ is Gaussian, and by integrating~\eqref{eq:6} it can be shown that its limiting variance agrees with the expression in~\cite{mbw-trim}.

\subsection{Max-sliced Wasserstein distance} \label{appwpp}
Unlike integration, the supremum functional is not smooth and does not possess a linear Hadamard derivative.
Write $MSW_{p}^p(P, Q) = \sup_{u \in \SS^{d-1}} W_p^p(P_u, Q_u)$, and note that $MSW_p^p(P, Q) = \omega(W)$, where $\omega:\ell^{\infty}(\SS^{d-1}) \to \RR$ is the supremum functional.
	It is shown in Theorem 2.1 of~\cite{crc-diff} that $\omega$ is Hadamard directionally differentiable with the derivative $$\omega'_f(g) = \mathop{\lim}\limits_{\epsilon \downarrow 0}\mathop{\sup}\limits_{x \in A_{\epsilon}(f)}g(x),$$ where $A_{\epsilon}(f) := \{x: f(x) \geq \sup f - \epsilon\}$. Moreover, if $f$ and $g$ are continuous on $\SS^{d-1}$ with respect to the standard Euclidean distance, then $\mathop{\lim}\limits_{\epsilon \downarrow 0} \mathop{\sup}\limits_{A_{\epsilon}(f)}g(x) = \mathop{\sup}\limits_{x \in A_0(f)}g(x)$ where $A_0(f) = \{x: f(x) = \sup f\}$. (See Corollary 2.3 of~\cite{crc-diff})

Applying the functional delta method to $\omega$ and the uniform weak convergence (\ref{eq:5}), we note that the limiting process $\GG$ has continuous samples paths a.s. so the limiting distribution of WPP can be written as \begin{equation} \label{eq:31}
    \omega'_{W_p^p(W)}(\GG) = \mathop{\sup}\limits_{u \in A_0(W_p^p(W))}\GG(u).
\end{equation}
In the spiked transport model (STM), this expression can be further simplified. The STM was introduced by~\cite{nr-stm} to formalize the situation where two distributions differ only in a low dimensional subspace of $\RR^d$. We describe the special case of one-dimensional spike here fix some $v \in \SS^{d-1}$ and let $X,\,Y \in L := \text{span}(v)$ be two random variables. Let $Z$ be another random variable independent of $(X,Y)$ and supported on the orthogonal complement $L^{\perp}$ of $L$. Then we define two distributions in $\RR^d$ by $P := \text{law}(X + Z),\,Q := \text{law}(Y + Z)$. \cite{nr-stm} show that $MSW_p^p(P,Q) = W_p(\text{law}(X),\text{law}(Y)) = W_2(P,Q)$. In addition, in this model the set $A_{\epsilon}(W_p^p(P_{\cdot},Q_{\cdot}))$ shrinks to the singleton set $\{v\}$ as $\epsilon$ goes down to $0$. In fact, the Hadamard derivative can be reduced to the random variable $\GG(u)$, i.e. the marginal of limiting Gaussian process $\GG$ along $v$. We summarize the result in the following theorem. 
\begin{theorem} \label{thm:4}
    Suppose that two compactly supported probability distributions $P$ and $Q$ in $\RR^d$ fit the spiked transport model with spike $v \in \SS^{d-1}$. Assume furthermore $P$ and $Q$ also satisfy \ref{CC}. Then, \begin{equation} \label{eq:32}
        \sqrt{n}\left(MSW_p^p(P_n,Q_n) - MSW_p^p(P,Q)\right) \stackrel{\text{d}}{\to} \GG(v).
    \end{equation}
\end{theorem}
\begin{remark} \label{rmk:4}
    Note that it is not necessary for $P$ and $Q$ to satisfy the spiked transport model in order to deduce the CLT of max-sliced Wasserstein. Namely, even if the set $\{u \in \SS^{d-1}: W_p^p(P_u,Q_u) = MSW(P,Q)\}$ is not a singleton, the same proof still works but the limiting distribution is the supremum of the Gaussian process $\GG$ over the set $\{u \in \SS^{d-1}: W_p^p(P_u,Q_u) = MSW(P,Q)\}$ which is not necessarily Gaussian. This example shows that certain functionals give rise to non-Gaussian limits, even though the limit in~\eqref{eq:5} is Gaussian.
\end{remark}
    
\section{Simulation Studies} \label{simulation}
\noindent We illustrate our distributional limit results in Monte Carlo simulations. Specifically, we investigate the speed of convergence for the empirical sliced Wasserstein distance and the empirical max-robust Wasserstein distance (also known as Wasserstein project pursuit~\cite{nr-stm})to their limit distributions. We also illustrate the accuracy of the approximation using the re-scaled bootstrap. All simulations were performed using Python. The Wasserstein distances as well as the sliced Wasserstein distances were calculated using the Python package $\textit{POT}$~\cite{flamary-pot} and the max-sliced Wasserstein distances were approximated by the Riemannian optimization method proposed in~\cite{lfh-prw}. \footnote{Code for simulations is available at \url{https://github.com/HelenXi/code_cltwd}.}
\subsection{Sliced Wasserstein Distance} \label{simusw}
\noindent We present an example that concerns two different distributions with connected projections along all directions. \\
\noindent Consider a simple model of transport where source and target distributions $P,\,Q$ are uniform on unit sphere $\SS^2$ and the unit sphere $\SS^2_{(1,1,1)}$ centered at $(1,1,1)$ respectively. \\
We first give an explicit representation of the theoretical limit of the example given in section \ref{simusw}. Fix any point $\theta \in \SS^{2}$, the projections of $P$ and $Q$ along $\theta$ are uniform over $(-1,1)$ and $(-1+a_{\theta},1+a_{\theta})$ respectively where $a_{\theta} := \theta_1 + \theta_2 + \theta_3$. Then the unique Kantorovich potential that achieves $2$-Wasserstein distance between $P_\theta$ and $Q_\theta$ is $\phi_{0}^{\theta}(x) = -2a_{\theta}x$. Hence, we have \begin{equation*}
    \sqrt{n}\left(W_2^2(P_{n\cdot},Q_{n\cdot} - W_2^2(P_{\cdot},Q_{\cdot})\right) \rightsquigarrow \GG,
\end{equation*} where $\GG$ is the mean-zero Gaussian process indexed by $\SS^2$ with covariance functions \begin{equation*}
    \EE\GG(u)\GG(v) = \frac{8}{3}a_ua_v\langle u,v \rangle.
\end{equation*}
It follows from Theorem \ref{thm:3} that the limiting distribution of the empirical $2$-Wasserstein distance is the centered Gaussian $S$ with variance \begin{equation*}
    \text{Var}(S) = \frac{8}{3} \int_{\SS^2}\int_{\SS^2}a_ua_v\langle u,v \rangle\,d\sigma(u)d\sigma(v) \approx 0.832.
\end{equation*}
\noindent We sample i.i.d. observations $X_1,\dots,X_n \sim P$ and $Y_1,\dots,Y_n \sim Q$ with size $n = 50,100,500$. This process is repeated $500$ times. We then compare the finite distributions of $1$-Wasserstein distance with the theoretical limit given in section \ref{appsw}. We demonstrate the results using kernel density estimators in Figure \ref{fig:1} along with the corresponding Q-Q plots. We see that the finite-sample empirical distribution gets closer to the limiting Gaussian distribution in \ref{eq:29} as the sample size $n$ increases. \\
\noindent In addition, we simulate the re-scaled plug-in bootstrap approximations by sampling $n = 1000$ observations of $P$ and $Q$. Fix some empirical SW $\sqrt{n}SW_2^2(P_n,Q_n)$, we generate $B = 500$ replications of $\sqrt{l}(SW_2^2(\hat{P}_n^*,\hat{Q}_n^*) - SW_2^2(P_n,Q_n))$. The distributions of the replications with various replacement numbers $l$, compared with the finite-sample empirical distribution and the theoretical limit, are shown in Figure \ref{fig:2}. We observe that the naive bootstrap ($l = n$) better approximates the finite sample distribution compared to fewer replacements ($l = n^{1/2},n^{3/4}$). This is consistent with the observation of inference on finite spaces.~\cite{oi-infproj} \\
\begin{figure} [!ht]
    \centering
    \begin{subfigure}{0.35\textwidth}
        \centering
        \includegraphics[width=1.\textwidth]{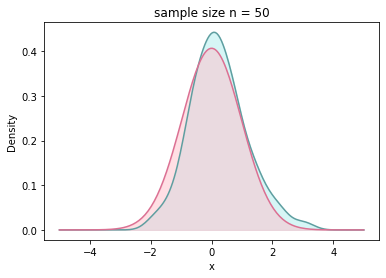} 
    \end{subfigure}
    \begin{subfigure}{0.35\textwidth}
        \centering
        \includegraphics[width=1.\textwidth]{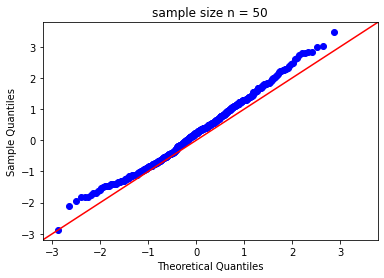} 
    \end{subfigure}
    \begin{subfigure}{0.35\textwidth}
        \centering
        \includegraphics[width=1.\textwidth]{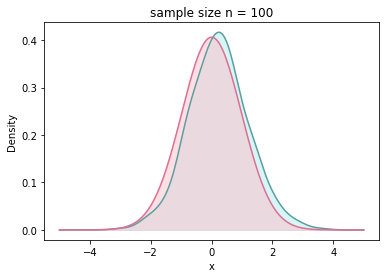} 
    \end{subfigure}
    \begin{subfigure}{0.35\textwidth}
        \centering
        \includegraphics[width=1.\textwidth]{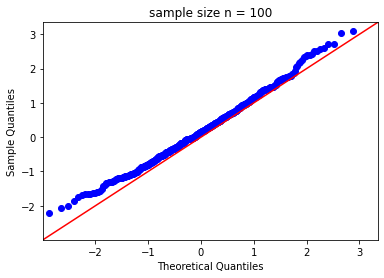} 
    \end{subfigure}
    \begin{subfigure}{0.35\textwidth}
        \centering
        \includegraphics[width=1.\textwidth]{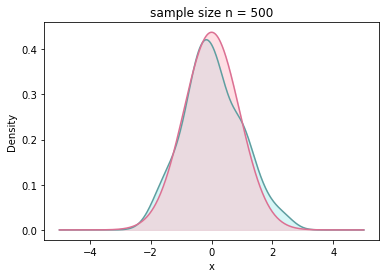} 
    \end{subfigure}
    \begin{subfigure}{0.35\textwidth}
        \centering
        \includegraphics[width=1.\textwidth]{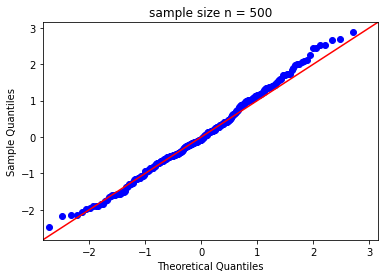} 
    \end{subfigure}
    \caption{Left: Comparison of the finite sample density (pale turquoise) and the limit distribution of the empirical sliced distance (pink). \\ Right: The corresponding Q-Q plots where the red solid line indicates perfect fit.}
    \label{fig:1}
\end{figure}
\begin{figure}[!ht]
    \centering
    \begin{subfigure}{0.33\textwidth}
        \centering
        \includegraphics[width=1.\textwidth]{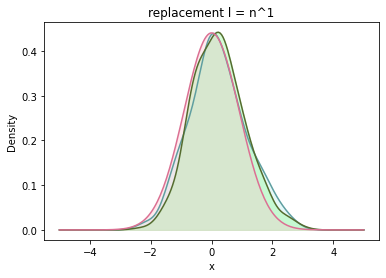} 
    \end{subfigure}
    \begin{subfigure}{0.33\textwidth}
        \centering
        \includegraphics[width=1.\textwidth]{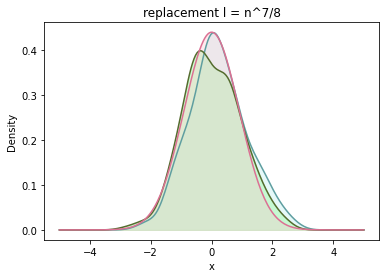} 
    \end{subfigure}
    \begin{subfigure}{0.33\textwidth}
        \centering
        \includegraphics[width=1.\textwidth]{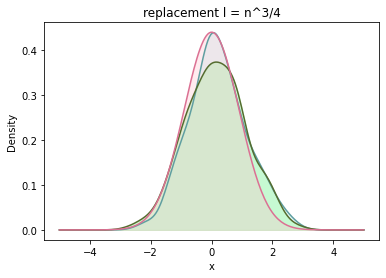} 
    \end{subfigure}
    \begin{subfigure}{0.33\textwidth}
        \centering
        \includegraphics[width=1.\textwidth]{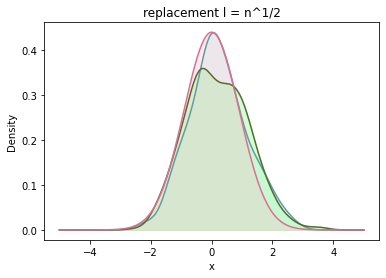} 
    \end{subfigure}
    \caption{Bootstrap for the empirical sliced distance. Illustration of the re-scaled plug-in bootstrap approximation ($n = 1000$) with replacement $l \in \{n, n^{7/8}, n^{4/5}, n^{1/2}\}$. Finite bootstrap densities (pale green) are compared to the corresponding finite sample density (pale turquoise) and the limit distribution (pink).}
    \label{fig:2}
\end{figure}

\subsection{Max-sliced Wasserstein distance}\label{simuwpp} 
\noindent We present an example that simulate different behavior of empirical Wasserstein projection pursuit when $p = 2$. \\
\noindent  We take $P$ to be the uniform distribution on the unit sphere $\SS^2$ and $Q$ to be uniform on the surface of ellipsoid $x^2/a^2 + y^2 + z^2 = 1$ where $a = 8.5$. We sample i.i.d. observations with size $n = 50,100,500$ and this process is repeated $2000$ times. The estimation plotted in the top part of Figure \ref{fig:3} indicates that the finite sample distributions approximate the limiting Gaussian distribution derived in Theorem \ref{thm:4} very well even when the sample size is small.\\
\noindent In terms of the re-scaled bootstrap, the accuracy of the bootstrap approximation seems to be good for all replacement numbers in this case, which is again consistent with the observation of the case when the underlying distributions are supported in finite sets.~\cite{oi-infproj} See Figure \ref{fig:4} for the simulation. \\
\begin{figure}[!ht]
    \centering
    \begin{subfigure}{0.3\textwidth}
        \centering
        \includegraphics[width=1.\textwidth]{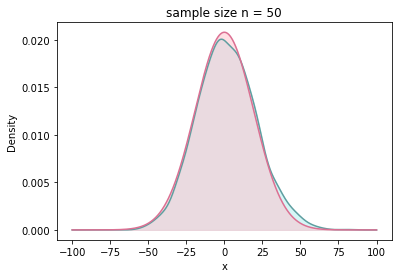}
    \end{subfigure}
    \begin{subfigure}{0.3\textwidth}
        \centering
        \includegraphics[width=1.\textwidth]{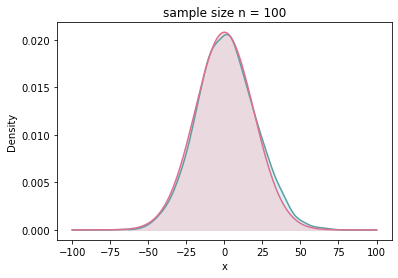} 
    \end{subfigure}
    \begin{subfigure}{0.3\textwidth}
        \centering
        \includegraphics[width=1.\textwidth]{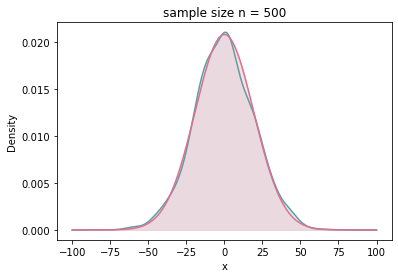} 
    \end{subfigure}
    \caption{Comparison of the finite sample density (pale turquoise) and the limit distribution of the empirical one-dimensional WPP (pink).}
    \label{fig:3}
\end{figure}
\begin{figure}[!ht]
    \centering
    \begin{subfigure}{0.33\textwidth}
        \centering
        \includegraphics[width=1.\textwidth]{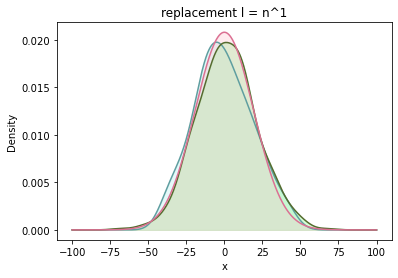} 
    \end{subfigure}
    \begin{subfigure}{0.33\textwidth}
        \centering
        \includegraphics[width=1.\textwidth]{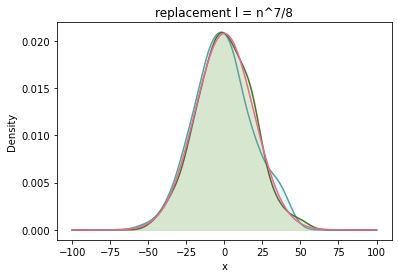} 
    \end{subfigure}
    \begin{subfigure}{0.33\textwidth}
        \centering
        \includegraphics[width=1.\textwidth]{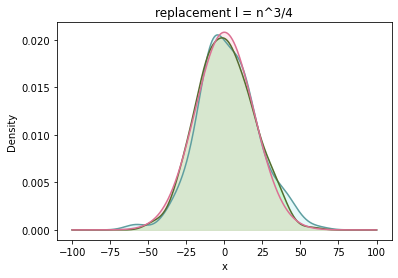} 
    \end{subfigure}
    \begin{subfigure}{0.33\textwidth}
        \centering
        \includegraphics[width=1.\textwidth]{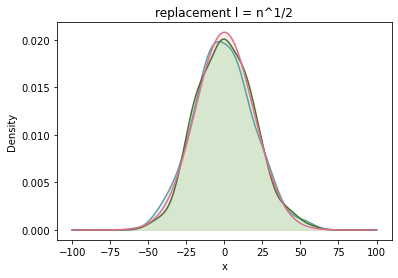} 
    \end{subfigure}
    \caption{Bootstrap for the empirical one-dimensional WPP. Illustration of the re-scaled plug-in bootstrap approximation ($n = 1000$) with replacement $l \in \{n, n^{7/8}, n^{4/5}, n^{1/2}\}$. Finite bootstrap densities (pale green) are compared to the corresponding finite sample density (pale turquoise) and the limit distribution (pink).}
    \label{fig:4}
\end{figure}

\clearpage

\bibliographystyle{alpha}
\bibliography{main}

\newcommand{\etalchar}[1]{$^{#1}$}
\begin{thebibliography}{MBNWW21}

\bibitem[AWR17]{AltWeeRig17}
Jason Altschuler, Jonathan Weed, and Philippe Rigollet.
\newblock Near-linear time approximation algorithms for optimal transport via
  {S}inkhorn iteration.
\newblock In {\em Advances in Neural Information Processing Systems 30: Annual
  Conference on Neural Information Processing Systems 2017, 4-9 December 2017,
  Long Beach, CA, {USA}}, pages 1961--1971, 2017.

\bibitem[BGT{\etalchar{+}}17]{bousquet2017optimal}
Olivier Bousquet, Sylvain Gelly, Ilya Tolstikhin, Carl-Johann Simon-Gabriel,
  and Bernhard Schoelkopf.
\newblock From optimal transport to generative modeling: the vegan cookbook.
\newblock {\em arXiv preprint arXiv:1705.07642}, 2017.

\bibitem[BL19]{bl-1d}
Sergey~G. Bobkov and Michel Ledoux.
\newblock One-dimensional empirical measures, order statistics, and kantorovich
  transport distances.
\newblock {\em Memoirs of the American Mathematical Society}, 2019.

\bibitem[BRPP15]{BonRabPey15}
Nicolas Bonneel, Julien Rabin, Gabriel Peyr\'{e}, and Hanspeter Pfister.
\newblock Sliced and {R}adon {W}asserstein barycenters of measures.
\newblock {\em J. Math. Imaging Vision}, 51(1):22--45, 2015.

\bibitem[CFTR17]{CouFlaTui17}
Nicolas Courty, R{\'{e}}mi Flamary, Devis Tuia, and Alain Rakotomamonjy.
\newblock Optimal transport for domain adaptation.
\newblock {\em {IEEE} Trans. Pattern Anal. Mach. Intell.}, 39(9):1853--1865,
  2017.

\bibitem[CRC20]{crc-diff}
Javier C'arcamo, Luis-Alberto Rodr'iguez, and Antonio Cuevas.
\newblock Directional differentiability for supremum-type functionals:
  Statistical applications.
\newblock {\em Bernoulli}, 2020.

\bibitem[Cut13]{Cut13}
Marco Cuturi.
\newblock Sinkhorn distances: Lightspeed computation of optimal transport.
\newblock In {\em Advances in Neural Information Processing Systems 26: 27th
  Annual Conference on Neural Information Processing Systems 2013. Proceedings
  of a meeting held December 5-8, 2013, Lake Tahoe, Nevada, United States.},
  pages 2292--2300, 2013.

\bibitem[DBGM99]{del1999central}
Eustasio Del~Barrio, Evarist Gin{\'e}, and Carlos Matr{\'a}n.
\newblock Central limit theorems for the wasserstein distance between the
  empirical and the true distributions.
\newblock {\em Annals of Probability}, pages 1009--1071, 1999.

\bibitem[dBGSL21]{del2021central}
Eustasio del Barrio, Alberto Gonz{\'a}lez-Sanz, and Jean-Michel Loubes.
\newblock Central limit theorems for general transportation costs.
\newblock {\em arXiv preprint arXiv:2102.06379}, 2021.

\bibitem[dBGSL22]{del2022central}
Eustasio del Barrio, Alberto Gonz{\'a}lez-Sanz, and Jean-Michel Loubes.
\newblock Central limit theorems for semidiscrete wasserstein distances.
\newblock {\em arXiv preprint arXiv:2202.06380}, 2022.

\bibitem[DBGU05]{del2005asymptotics}
Eustasio Del~Barrio, Evarist Gin{\'e}, and Frederic Utzet.
\newblock Asymptotics for l2 functionals of the empirical quantile process,
  with applications to tests of fit based on weighted wasserstein distances.
\newblock {\em Bernoulli}, 11(1):131--189, 2005.

\bibitem[DBL19]{bl-clt}
Eustasio Del~Barrio and Jean-Michel Loubes.
\newblock Central limit theorems for empirical transportation cost in general
  dimension.
\newblock {\em The Annals of Probability}, 47(2):926--951, 2019.

\bibitem[DGS21]{DebGhoSen21}
Nabarun Deb, Promit Ghosal, and Bodhisattva Sen.
\newblock Rates of estimation of optimal transport maps using plug-in
  estimators via barycentric projections.
\newblock {\em Advances in Neural Information Processing Systems}, 34, 2021.

\bibitem[DHS{\etalchar{+}}19]{deshpande2019max}
Ishan Deshpande, Yuan-Ting Hu, Ruoyu Sun, Ayis Pyrros, Nasir Siddiqui, Sanmi
  Koyejo, Zhizhen Zhao, David Forsyth, and Alexander~G Schwing.
\newblock Max-sliced wasserstein distance and its use for gans.
\newblock In {\em Proceedings of the IEEE/CVF Conference on Computer Vision and
  Pattern Recognition}, pages 10648--10656, 2019.

\bibitem[DM09]{dm-cover}
Oleksiy Dovgoshey and Olli Martio.
\newblock Products of metric spaces, covering numbers, packing numbers and
  characterizations of ultrametric spaces.
\newblock {\em arXiv: Metric Geometry}, 2009.

\bibitem[Dud69]{Dud69}
R.~M. Dudley.
\newblock The speed of mean {G}livenko-{C}antelli convergence.
\newblock {\em Ann. Math. Statist}, 40:40--50, 1969.

\bibitem[D{\"u}m93]{dumbgen1993nondifferentiable}
Lutz D{\"u}mbgen.
\newblock On nondifferentiable functions and the bootstrap.
\newblock {\em Probability Theory and Related Fields}, 95(1):125--140, 1993.

\bibitem[FCG{\etalchar{+}}21]{flamary-pot}
R{\'e}mi Flamary, Nicolas Courty, Alexandre Gramfort, Mokhtar~Z. Alaya,
  Aur{\'e}lie Boisbunon, Stanislas Chambon, Laetitia Chapel, Adrien Corenflos,
  Kilian Fatras, Nemo Fournier, L{\'e}o Gautheron, Nathalie~T.H. Gayraud,
  Hicham Janati, Alain Rakotomamonjy, Ievgen Redko, Antoine Rolet, Antony
  Schutz, Vivien Seguy, Danica~J. Sutherland, Romain Tavenard, Alexander Tong,
  and Titouan Vayer.
\newblock Pot: Python optimal transport.
\newblock {\em Journal of Machine Learning Research}, 22(78):1--8, 2021.

\bibitem[GKRS22]{goldfeld2022statistical}
Ziv Goldfeld, Kengo Kato, Gabriel Rioux, and Ritwik Sadhu.
\newblock Statistical inference with regularized optimal transport.
\newblock {\em arXiv preprint arXiv:2205.04283}, 2022.

\bibitem[GPC18]{GenPeyCut18}
Aude Genevay, Gabriel Peyr{\'{e}}, and Marco Cuturi.
\newblock Learning generative models with sinkhorn divergences.
\newblock In {\em International Conference on Artificial Intelligence and
  Statistics, {AISTATS} 2018, 9-11 April 2018, Playa Blanca, Lanzarote, Canary
  Islands, Spain}, pages 1608--1617, 2018.

\bibitem[GS19]{GhoSen19}
Promit Ghosal and Bodhisattva Sen.
\newblock Multivariate ranks and quantiles using optimal transport:
  Consistency, rates, and nonparametric testing.
\newblock {\em arXiv preprint arXiv:1905.05340}, 2019.

\bibitem[HKSM22]{hundrieser2022unifying}
Shayan Hundrieser, Marcel Klatt, Thomas Staudt, and Axel Munk.
\newblock A unifying approach to distributional limits for empirical optimal
  transport.
\newblock {\em arXiv preprint arXiv:2202.12790}, 2022.

\bibitem[LFH{\etalchar{+}}20]{lfh-prw}
Tianyi Lin, Chenyou Fan, Nhat Ho, Marco Cuturi, and Michael Jordan.
\newblock Projection robust wasserstein distance and riemannian optimization.
\newblock {\em Advances in Neural Information Processing Systems},
  33:9383--9397, 2020.

\bibitem[MBNWW21]{manole2021plugin}
Tudor Manole, Sivaraman Balakrishnan, Jonathan Niles-Weed, and Larry Wasserman.
\newblock Plugin estimation of smooth optimal transport maps.
\newblock {\em arXiv preprint arXiv:2107.12364}, 2021.

\bibitem[MBW19]{mbw-trim}
Tudor Manole, Sivaraman Balakrishnan, and Larry Wasserman.
\newblock Minimax confidence intervals for the sliced wasserstein distance.
\newblock {\em arXiv: Statistics Theory}, 2019.

\bibitem[MNW21]{manole2021sharp}
Tudor Manole and Jonathan Niles-Weed.
\newblock Sharp convergence rates for empirical optimal transport with smooth
  costs.
\newblock {\em arXiv preprint arXiv:2106.13181}, 2021.

\bibitem[NDC{\etalchar{+}}20]{nadjahi2020statistical}
Kimia Nadjahi, Alain Durmus, L{\'e}na{\"\i}c Chizat, Soheil Kolouri, Shahin
  Shahrampour, and Umut Simsekli.
\newblock Statistical and topological properties of sliced probability
  divergences.
\newblock {\em Advances in Neural Information Processing Systems},
  33:20802--20812, 2020.

\bibitem[NHPB20]{nguyen2020distributional}
Khai Nguyen, Nhat Ho, Tung Pham, and Hung Bui.
\newblock Distributional sliced-wasserstein and applications to generative
  modeling.
\newblock {\em arXiv preprint arXiv:2002.07367}, 2020.

\bibitem[NWR19]{nr-stm}
Jonathan Niles-Weed and Philippe Rigollet.
\newblock Estimation of wasserstein distances in the spiked transport model.
\newblock {\em arXiv preprint arXiv:1909.07513}, 2019.

\bibitem[OI22a]{okano2022inference}
Ryo Okano and Masaaki Imaizumi.
\newblock Inference for projection-based wasserstein distances on finite
  spaces.
\newblock {\em arXiv preprint arXiv:2202.05495}, 2022.

\bibitem[OI22b]{oi-infproj}
Ryotaro Okano and Masaaki Imaizumi.
\newblock Inference for projection-based wasserstein distances on finite
  spaces.
\newblock 2022.

\bibitem[PC19]{paty2019subspace}
Fran{\c{c}}ois-Pierre Paty and Marco Cuturi.
\newblock Subspace robust wasserstein distances.
\newblock In {\em International conference on machine learning}, pages
  5072--5081. PMLR, 2019.

\bibitem[RHS17]{RedHabSeb17}
Ievgen Redko, Amaury Habrard, and Marc Sebban.
\newblock Theoretical analysis of domain adaptation with optimal transport.
\newblock In {\em Joint European Conference on Machine Learning and Knowledge
  Discovery in Databases}, pages 737--753. Springer, 2017.

\bibitem[R{\"o}m06]{rom-delta}
Werner R{\"o}misch.
\newblock {\em Delta Method, Infinite Dimensional}.
\newblock John Wiley $\&$ Sons, Ltd, 2006.

\bibitem[RPDB11]{rpdb-sliced}
Juline Rabin, Gabriel Peyr{\'e}, Julie Delon, and Marc Bernot.
\newblock Wasserstein barycenter and its application to texture mixing.
\newblock {\em International Conference on Scale Space and Variational Methods
  in Computer Vision}, pages 435--446, 2011.

\bibitem[Rud86]{rud-analysis}
Walter Rudin.
\newblock {\em Principles of Mathematical Analysis}.
\newblock McGraw - Hill Book C., 1986.

\bibitem[San15]{San15}
Filippo Santambrogio.
\newblock {\em Optimal transport for applied mathematicians}, volume~87 of {\em
  Progress in Nonlinear Differential Equations and their Applications}.
\newblock Birkh\"{a}user/Springer, Cham, 2015.
\newblock Calculus of variations, PDEs, and modeling.

\bibitem[Sha90]{sha-gatu}
Arnold~S. Shapiro.
\newblock On concepts of directional differentiability.
\newblock {\em Journal of optimization theory and applications},
  66(3):477--487, 1990.

\bibitem[SHM22]{sha-kanto}
Thomas Staudt, Shayan Hundrieser, and Axel Munk.
\newblock On the uniqueness of kantorovich potentials.
\newblock 2022.

\bibitem[SM18]{sommerfeld2018inference}
Max Sommerfeld and Axel Munk.
\newblock Inference for empirical wasserstein distances on finite spaces.
\newblock {\em Journal of the Royal Statistical Society: Series B (Statistical
  Methodology)}, 80(1):219--238, 2018.

\bibitem[SP18]{SinPoc18}
Shashank Singh and Barnab{\'a}s P{\'o}czos.
\newblock Minimax distribution estimation in {W}asserstein distance.
\newblock {\em arXiv preprint arXiv:1802.08855}, 2018.

\bibitem[SST{\etalchar{+}}19]{schiebinger2019optimal}
Geoffrey Schiebinger, Jian Shu, Marcin Tabaka, Brian Cleary, Vidya Subramanian,
  Aryeh Solomon, Joshua Gould, Siyan Liu, Stacie Lin, Peter Berube, et~al.
\newblock Optimal-transport analysis of single-cell gene expression identifies
  developmental trajectories in reprogramming.
\newblock {\em Cell}, 176(4):928--943, 2019.

\bibitem[TSM19]{tameling2019empirical}
Carla Tameling, Max Sommerfeld, and Axel Munk.
\newblock Empirical optimal transport on countable metric spaces:
  Distributional limits and statistical applications.
\newblock {\em The Annals of Applied Probability}, 29(5):2744--2781, 2019.

\bibitem[vdVW96]{vw-ep}
Aad van~der Vaart and Jon~August Wellner.
\newblock Weak convergence and empirical processes: With applications to
  statistics.
\newblock 1996.

\bibitem[Vil08]{vil-ot}
C.~Villani.
\newblock {\em Optimal Transport: Old and New}.
\newblock Grundlehren der mathematischen Wissenschaften. Springer Berlin
  Heidelberg, 2008.

\bibitem[WB19]{weedBach}
Jonathan Weed and Francis Bach.
\newblock {Sharp asymptotic and finite-sample rates of convergence of empirical
  measures in Wasserstein distance}.
\newblock {\em Bernoulli}, 25(4A):2620 -- 2648, 2019.

\bibitem[Wel05]{wellner-emp}
Jon~A. Wellner.
\newblock Empirical processes: Theory and applications.
\newblock 2005.

\bibitem[XH22]{XuHua22}
Xianliang Xu and Zhongyi Huang.
\newblock Central limit theorem for the sliced 1-wasserstein distance and the
  max-sliced 1-wasserstein distance.
\newblock 05 2022.

\bibitem[YDV{\etalchar{+}}20]{dai2018autoencoder}
Karren~Dai Yang, Karthik Damodaran, Saradha Venkatachalapathy, Ali~C
  Soylemezoglu, GV~Shivashankar, and Caroline Uhler.
\newblock Predicting cell lineages using autoencoders and optimal transport.
\newblock {\em PLoS computational biology}, 16(4):e1007828, 2020.

\bibitem[YNN{\etalchar{+}}21]{yang2021optimal}
Yunan Yang, Levon Nurbekyan, Elisa Negrini, Robert Martin, and Mirjeta Pasha.
\newblock Optimal transport for parameter identification of chaotic dynamics
  via invariant measures.
\newblock {\em arXiv preprint arXiv:2104.15138}, 2021.

\end{thebibliography}

\clearpage

\appendix
\section{Proofs}
\subsection{Proof of Theorem \ref{thm:1}}
To prove \cref{thm:1}, we will first show a uniform central limit theorem for the empirical process indexed by the set $\sC \times \SS^{d-1}$ where $\sC := \{f: [-R,R] \to \RR,\,\phi(0) = 0,\,\|f\|_{\text{Lip}} \leq L\}$ for some positive quantity $L$.
Explicitly, for $f \in \sC$, $u \in \SS^{d-1}$, write $h_u$ for the function $x \mapsto u^\top x$ and define $\BB_{n} \in \ell^\infty(\sC \times \SS^{d-1})$ by
\begin{align*}
	\BB_{nm}(f, u) & := \sqrt n(P_n - P)(f \circ h_u) + \sqrt m(Q_m - Q)(f^c \circ h_u) \\
	& = \frac{1}{\sqrt n} \sum_{i=1}^n f(u^\top X_i) - \EE f(u^\top X_i) + \frac{1}{\sqrt m} \sum_{j=1}^m f^c(u^\top Y_i) - \EE f^c(u^\top Y_i)\,.
\end{align*}
We first show that $\BB_n$ possesses a weak limit.
\begin{proposition}\label{prop:weak}
	As $n, m \to \infty$,
	the empirical process $\BB_{nm}$ satisfies
	\begin{equation*}
		\BB_{nm} \rightsquigarrow \BB \quad \text{in $\ell^\infty(\sC \times \SS^{d-1})$,}
	\end{equation*}
	where $\BB$ is the tight Gaussian process with covariance
	\begin{multline}
		\EE \BB(f, u) \BB(g, v) = \int f(u^\top x) g(v^\top x) \dd P(x) - \int f(u^\top x) \dd P(x) \int g(v^\top x) \dd P(x) \\
		+  \int f^c(u^\top y) g^c(v^\top y) \dd Q(y) - \int f^c(u^\top y) \dd Q(y) \int g^c(v^\top y) \dd Q(y)\,.
	\end{multline}
	Moreover, this process is uniformly continuous with respect to the semimetric
	\begin{equation}\label{semi_metric}
		\rho((f, u), (g, v)) = \|f \circ h_u - g \circ h_v\|_{L^2(P)} + \|f^c \circ h_u - g^c \circ h_v\|_{L^2(Q)}\,,
	\end{equation}
	with respect to which $\sC \times \SS^{d-1}$ is totally bounded.
\end{proposition}
\begin{proof}
	The assertions of this proposition will follow from the fact that the classes of functions $\sF := \{f \circ h_u(x) : (f, u) \in \sC \times \SS^{d-1}\}$ and $\sF^C :=\{f^c \circ h_u(x) : (f, u) \in \sC \times \SS^{d-1}\}$ are $P$ and $Q$-Donsker, respectively.
	Indeed, if we assume this Donsker property, then we have
	\begin{align*}
		\sqrt{n}(P_n - P) & \rightsquigarrow \BB_P \quad \text{in $\ell^\infty(\sF)$} \\
		\sqrt{m}(Q_m - Q) & \rightsquigarrow \BB_Q \quad \text{in $\ell^\infty(\sF^C)$}
	\end{align*}
	for tight $P$- and $Q$-Brownian bridges $\BB_P$ and $\BB_Q$, respectively~\cite[see][Section 2.1]{vw-ep}.
	These tight Gaussian processes possess uniformly continuous sample paths with respect to the semi-metrics $\rho_P$ and $\rho_Q$, respectively, where for $F, G \in \sF$,
	\begin{equation*}
		\rho^2_P(F, G) = \int (F-G)^2 \dd P - \left(\int (F-G) \dd P\right )^2\,,
	\end{equation*}
	and analogously for $Q$, and $\sF$ and $\sF^C$ are totally bounded with respect to these semi-metrics~\cite[Example 1.5.10]{vw-ep}.
	In particular, since $\rho_P$ is dominated by the $L^2(P)$ norm, and likewise for $Q$, these processes are also uniformly $L^2(P)$ and $L^2(Q)$ continuous and $\sF$ and $\sF^C$ are $L^2(P)$ and $L^2(Q)$ totally bounded.
	
	Since $P_n$ and $Q_m$ are independent, the above considerations imply that
	\begin{equation*}
		(\sqrt{n}(P_n - P), \sqrt{m}(Q_m - Q)) \rightsquigarrow (\BB_P, \BB_Q) \quad \text{in $\ell^\infty(\sF) \times \ell^\infty(\sF^C)$,}
	\end{equation*}
	for a tight Gaussian limit $(\BB_P, \BB_Q)$ with sample paths almost surely continuous with respect to the metric on $\sF \times \sF^C$ given by the sum of the $L^2(P)$ and $L^2(Q)$ metrics on $\sF$ and $\sF^C$.
	Finally, there exists a continuous map from $\ell^\infty(\sF) \times \ell^\infty(\sF^C)$ to $\ell^\infty(\sC \times \SS^{d-1})$ given by associating $(S, T) \in \ell^\infty(\sF) \times \ell^\infty(\sF^C)$ with the element of $\ell^\infty(\sC \times \SS^{d-1})$ sending $(f, u)$ to $S(f \circ h_u) + T(f^c \circ h_u)$, and the continuous mapping theorem therefore furnishes the desired convergence.
	
	It remains to show that $\sF$ and $\sF^C$ are $P$- and $Q$-Donsker.
	We first prove that $\sF$ is $P$-Donsker.
	By \cite[Theorem 2.5.6]{vw-ep}, it suffices to show that
	\begin{equation}
		\int_0^{\infty} \sqrt{\log N_{[]}(\varepsilon,\sF,L_2(P))} \dd \epsilon < \infty\,.
	\end{equation}
	
	By \cref{lemma:2}, we may replace the bracketing number $N_{[]}(\varepsilon,\sF,L_2(P))$ by the uniform covering number $N(\varepsilon/2,\sF,\|\cdot\|_{\infty})$.
	
	\cite[Theorem 2.7.1]{vw-ep} gives an upper bound for the covering entropy of $\sC$: there exists some positive constant $C_1$ depending on $R$ and $L$ such that \begin{equation}
		\log N(\varepsilon,\sC,\|\cdot\|_{\infty}) \leq C_1 \varepsilon^{-1}.
	\end{equation}
	\cite[Lemma 6.2]{wellner-emp} shows that there exists a positive constant $C_2$ depending on $R$ such that \begin{equation}
		N(\varepsilon,\{h_u,\,u \in \SS^{d-1}\},\|\cdot\|_{\infty}) = N(\varepsilon/2,\SS^{d-1},\|\cdot\|_{2}) \leq C_2 \varepsilon^{-d}\,.
	\end{equation}
	Consequently, apply \cref{lemma:3} to $N(\varepsilon,\sF,\|\cdot\|_{\infty})$ and we get
	\begin{equation}\label{uniform_covering_bound}
		N(\varepsilon,\sF,\|\cdot\|_{\infty}) \leq \log N(\varepsilon,\sC,\|\cdot\|_{\infty}) + \log N(\varepsilon/2,\SS^{d-1},\|\cdot\|_{2}) \leq C_1 \varepsilon^{-1} + \log( C_2 \varepsilon^{-d})\,.
	\end{equation}
	We obtain that
	\begin{align*}
		\int_0^{\infty} \sqrt{\log N_{[]}(\varepsilon,\sF,L_2(P))} \dd \epsilon & \leq \int_0^{\infty} \sqrt{N(\varepsilon/2,\sF,\|\cdot\|_{\infty})} \dd \epsilon \\
		& = \int_0^{2\mathrm{diam}(\sF)}  \sqrt{N(\varepsilon/2,\sF,\|\cdot\|_{\infty})} \dd \epsilon \\
		& \leq \int_0^{2\mathrm{diam}(\sF)} \sqrt{C_1 \varepsilon^{-1} + \log( C_2 \varepsilon^{-d})} \dd \epsilon \\
		& < \infty\,. 
	\end{align*}
	Therefore $\sF$ is $P$-Donsker.
	The argument for $\sF^C$ is identical, since \cref{lemma:4} shows that the estimate~\eqref{uniform_covering_bound} holds for $\sF^C$ as well.
\end{proof}

To prove \cref{thm:1}, we combine the above result with the functional delta method.
Let $\iota:\ell^{\infty}(\sC \times \SS^{d-1}) \to \ell^{\infty}(\SS^{d-1})$ be defined by 
\begin{equation}\label{iota_def}
	\iota(\Phi)(u) := \sup_{f \in \sC}\Phi(f, u)\,.
\end{equation}
The following proposition shows that $\iota$ is Hadamard directionally differentiable tangentially to the set of continuous functions at all functions for which the supremum in~\eqref{iota_def} is uniquely achieved.

\begin{proposition}\label{prop:hadamard}
	Let $\iota: \ell^{\infty}(\sC \times \SS^{d-1}) \to \ell^{\infty}(\SS^{d-1})$ be defined as in \eqref{iota_def}, and denote by $\mathsf C_u(\sC \times \SS^{d-1}, \rho)$ the set of elements of $\ell^{\infty}(\sC \times \SS^{d-1})$ which are uniformly continuous with respect to the semi-metric $\rho$ defined in~\eqref{semi_metric}.
	Then for $\Phi \in \mathsf C_u(\sC \times \SS^{d-1}, \rho)$ such that $\Phi(\cdot,u)$ has a unique maximizer for every $u \in \SS^{d-1}$, the function $\iota$ is Hadamard directionally differentiable at $\Phi$ tangentially to $\mathsf C_u(\sC \times \SS^{d-1}, \rho)$, with derivative $\iota'_{\Phi}: \ell^{\infty}(\sC \times \SS^{d-1}) \to \ell^{\infty}(\SS^{d-1})$ given by
	\begin{equation}
		\iota'_{\Phi}(\Psi)(u) =\Psi(f_u, u)\,, \quad u \in \SS^{d-1}, \Psi \in\mathsf C_u(\sC \times \SS^{d-1}, \rho) \,,
	\end{equation}
	where $\Phi(f_u,u) = \sup_\sC \Phi(\cdot,u)$.
\end{proposition}
\begin{proof}
	Fix an arbitrary $u \in \SS^{d-1}$.
	\cite[Theorem 2.1]{crc-diff} shows that the function $\iota^u: \ell^\infty(\sC \times \SS^{d-1}) \to \RR$ defined by
	\begin{equation*}
		\iota^u(\Phi) := \sup_{f \in \sC} \Phi(f, u)
	\end{equation*}
	is Hadamard directionally differentiable, with derivative
	\begin{equation}\label{der}
		{(\iota^u)}_\Phi'(\Psi) = \lim_{\epsilon \to 0} \sup_{f \in B_{\epsilon, u}(\Phi)} \Psi(f, u) 
	\end{equation}
	where $B_{\epsilon, u}(\Phi) := \{f \in \sC: \Phi(f, u) \geq \sup_\sC \Phi(\cdot, u) - \epsilon\}$.
	Moreover, we claim that if $\Phi$ and $\Psi$ are uniformly continuous, 
	then the expression for the derivative simplifies to
	\begin{equation}\label{easy_der}
		{(\iota^u)}_\Phi'(\Psi) = \Psi(f_u, u) = \iota'_\Phi(\Psi)(u)\,.
	\end{equation}
	To see this, we define $\phi, \psi \in \ell^\infty(\sC)$ by $\phi(\cdot) = \Phi(\cdot, u)$ and $\psi(\cdot) = \Psi(\cdot, u)$, so that the right side of \eqref{der} reads
	\begin{equation*}
		\lim_{\epsilon \to 0} \sup_{f \in B_\epsilon(\phi)} \psi(f)\,,
	\end{equation*}
	where $B_\epsilon(\phi) := \{f \in \sC: \phi(f) \geq \sup_{\sC} \phi - \epsilon\}$.
	By assumption, the functions $\phi$ and $\psi$ are uniformly continuous with respect to the semi-metric
		\begin{equation*}
			\rho^u(f, g) := \|f \circ h^u - g \circ h^u\|_{L^2(P)} + \|f^c \circ h^u - g^c \circ h^u\|_{L^2(Q)}\,,
		\end{equation*}
 		and $\sC$ is totally bounded with respect to this semi-metric.
 		Following \cite[Corollary 2.5]{crc-diff}, it is enough to show that $\sC$ is complete with respect to $\rho^u$.
 		The completeness of $L^2((h_u)_\sharp P)$ implies that any sequence $f_n \in \sC$ which is Cauchy with respect to $\rho^u$ possesses a limit $f$, and by passing to a subsequence we may assume that $f_n \to f$ pointwise, and, since the elements of $\sC$ are bounded and equicontinuous, we may further assume that $f_n \to f$ uniformly on $[-R, R]$ by the Arzel\`a--Ascoli theorem.
 		Since $\sC$ is closed with respect to pointwise convergence, $f \in \sC$, and since $c$-transforms are preserved under uniform convergence, we also have $f_{n_k}^c \to f^c$ uniformly in $[-R, R]$.
 		Therefore $(f_{n_k}, f_{n_k}^c) \to (f, f^c)$ for some $f \in \sC$ uniformly, and hence $f_{n_k} \to f$ in $\rho^u$.
		This proves~\eqref{easy_der}.
	
	We now turn to the differentiability of $\iota$.
	To prove the Hadamard differentiability of $\iota$, according to Proposition 3.5 of~\cite{sha-gatu}, is equivalent to prove $\iota$ is Lipschitz and that the function $\iota_n(\Phi, \Psi) \in \ell^\infty(\SS^{d-1})$ defined by
	\begin{equation*}
		\iota_n(\Phi,\Psi)(\cdot) := \mathop{\sup}\limits_{f \in \sC}(s_n\Phi(f,\cdot) + \Psi(f,\cdot)) - s_n\mathop{\sup}\limits_{f\in \sC}\Phi(f,\cdot)
	\end{equation*}
	converges uniformly to the limit $\iota'_\Phi(\Psi)(\cdot)$ for any positive increasing sequence $s_n \to \infty$. The Lipschitz property is obvious. Indeed, for any $\Phi_1,\Phi_2 \in \ell^{\infty}(\sC \times \SS^{d-1})$,
	\begin{equation*}
		\mathop{\sup}\limits_{u \in \SS^{d-1}}\left|\mathop{\sup}\limits_{f}\Phi_1(f, u) - \mathop{\sup}\limits_{f}\Phi_2({f,u})\right| \leq \mathop{\sup}\limits_{f, u}\left|\Phi_1(f, u) - \Phi_2(f, u)\right| \leq \|\Phi_1 - \Phi_2\|_{\ell^{\infty}(\sC \times \SS^{d-1]})}.
	\end{equation*}
	
	As for the uniform convergence, we first show that $\iota_n(\Phi, \Psi)(u) \to \iota'_\Phi(\Psi)(u)$ pointwise.
	This follows directly from the Hadamard differentiability of $\iota^u$, since
	\begin{equation*}
		\iota_n(\Phi, \Psi)(u) = \iota^u(s_n \Phi + \Psi) - s_n \iota^u(\Phi) \to {(\iota^u)}_\Phi'(\Psi) = \iota'_\Phi(\Psi)(u) \quad \text{as $n \to \infty$.}
	\end{equation*}
	Moreover, we show in \cref{lem:continuous} that the functions $\iota_n(\Phi, \Psi)$ and $\iota'(\Phi,\Psi)$ are continuous on $\SS^{d-1}$. 
	Therefore, by~\cite[Theorem 7.13]{rud-analysis}, to show uniform convergence on $\SS^{d-1}$ it suffices to show that  the sequence $\{\iota_n(\Phi,\Psi)\}_{n \geq 1}$ is monotonically non-increasing for all $u \in \SS^{d-1}$. This follows directly from the definition of $\iota_n$. Indeed, for any $u \in \SS^{d-1}$, we have
	\begin{multline*}
		 \iota_n(\Phi,\Psi)(u) - \iota_{n+1}(\Phi,\Psi)(u) \\
		 = \mathop{\sup}\limits_{f \in \sC}\left(s_n\Phi(f, u) + \Psi(f, u)\right) - \mathop{\sup}\limits_{f \in \sC }\left(s_{n+1}\Phi(f, u) + \Psi(f, u)\right) + \mathop{\sup}\limits_{f \in \sC} (s_{n+1} - s_n)\Phi(f,u) \\
		 \geq \mathop{\sup}\limits_{f \in \sC}\left((s_n + s_{n+1} - s_n)\Phi(f,u) + \Psi(f,u)\right) - \mathop{\sup}\limits_{f \in \sC}\left(s_{n+1}\Phi(f, u) + \Psi(f, u)\right) = 0.
	\end{multline*} The last inequality results from the reverse triangle inequality for the supremum. This finishes the proof for Hadamard directional differentiability of $\iota$.
\end{proof}
\begin{lemma}\label{lem:continuous}
	For $\Phi, \Psi \in \mathsf C_u(\sC \times \SS^{d-1})$, then the functions $\iota_n(\Phi, \Psi)$ are continuous on $\SS^{d-1}$. If $\Phi(\cdot,u)$ has a unique maximizer for every $u \in \SS^{d-1}$, then $\iota'(\Phi,\Psi)$ is also continuous. 
\end{lemma}
\begin{proof}
	For the continuity of $\iota_n$, it suffices to show that if $u \to v$, then
	\begin{align}\label{phi_uc}
		\sup_{f \in \sC} |\Phi(f, u) - \Phi(f, v)| & \to 0\,,
	\end{align}
	and analogously for $\Psi$.
	This follows directly from uniform continuity: for any $\epsilon > 0$, there exist a $\delta > 0$ such that
	\begin{equation*}
		\rho((f, u), (g, v)) \leq \delta \implies |\Phi(f, u) - \Phi(g, v)| \leq \epsilon\,.
	\end{equation*}
	In particular, we have $\sup_{f \in \sC} |\Phi(f, u) - \Phi(f, v)|  \leq \epsilon$ if $\sup_f \rho((f, u), (f, v)) \leq \delta$.
	Moreover, since the elements of $\sC$ and their $c$-transforms are uniformly Lipschitz, we have
	\begin{equation*}
		\sup_{f \in \sC} \rho((f, u), (f, v)) \leq C \left(\|h_u - h_v\|_{L^2(P)} + \|h_u - h_v\|_{L^2(Q)}\right)
	\end{equation*}
	for a positive $C$ independent of $f$, and the right side of the above expression converges to $0$ as $u \to v$.
	Therefore~\eqref{phi_uc} holds, as does the analogous convergence for $\Psi$.
	This proves continuity of $\iota_n$.	
	
	For $\iota'$, we have \begin{align*}
		\iota_{\Phi}'(\Psi)(u) - \iota_{\Phi}(\Psi)(v) & = \Psi(f_u,u) - \Psi(f_v,v)
	\end{align*}
	where $f_u,f_v$ are the maximizers of $\Phi(\cdot,u),\Phi(\cdot,v)$ respectively.
	Choose any sequence $v_n \to v$ and let $f_n$ denote the unique maximizers of $\Phi(\cdot,v_n)$ correspondingly. 
	Since $\sC \times \SS^{d-1}$ is totally bounded and complete, we may upon passing to a subsequence assume $(f_n,v_n) \to (f,v) \in \sC \times \SS^{d-1}$, and by the uniform continuity of $\Phi$ we must have $\Phi(f,v) = \sup_{f \in \sC}\Phi(f,v)$, and since we have assumed that the supremum is uniquely achieved, $f = f_v$.
	Since this argument holds on any subsequence, we obtain that the whole sequence converges to $(f_v, v)$, and the uniform continuity of $\Psi$ implies that
 \begin{equation*}
		\lim_{n \to \infty} \Psi(f_n,v_n) = \Psi(f_v,v)\,,
	\end{equation*}
as desired.
\end{proof}
We are now in a position to prove the main theorem.
\begin{proof}[Proof of \cref{thm:1}]
	For the sake of notational simplicity, we prove the special case of (\ref{eq:5}) when $n = m$ with both sides multiplied by $\sqrt{2}$. Namely, we are going to show \begin{equation} \label{eq:12}
    	\sqrt{n}\left(W_p^p(P_{n\cdot},Q_{n\cdot}) - W_p^p(P_{\cdot},Q_{\cdot})\right) \rightsquigarrow \sqrt{2}\GG \quad \text{ in } \ell^{\infty}(\SS^{d-1}).
\end{equation}
	The general conclusion with $n \neq m$ follows by an analogous argument. 
	
	Fix $u \in \SS^{d-1}$.
	By Kantorovich duality, we may write the Wasserstein distance as
	\begin{equation*}
		W_p^p(P_u, Q_u) = \sup_{f \in \sC} \int f \circ h_u \dd P + \int f^c \circ h_u \dd Q\,.
	\end{equation*}
	Define $\Phi_{(P, Q)} : \sC \times \SS^{d-1} \to \RR$ by $\Phi_{(P,Q)} = \EE_{X \sim P}f \circ h_u(X) + \EE_{X \sim Q}f^c \circ h_u(Y)$. Note that $\Phi_{(P,Q)}$ is uniformly continuous with respect to $\rho$. Indeed, for any $f,g \in \sC$ and $u,v \in \SS^{d-1}$, \begin{align*}
		& \left|\Phi_{(P,Q)}(f,u) - \Phi_{(P,Q)}(g,v)\right|\\
		& \leq \EE_{X \sim P}\left|f \circ h_u(X) - g \circ h_v(X)\right| + \EE_{Y \sim Q}\left|f^c \circ h_u(Y) - g^c \circ h_v(Y)\right| \\
		& \leq \left(\EE_{X \sim P}\left|f \circ h_u(X) - g \circ h_v(X)\right|^2\right)^{1/2} + \left(\EE_{Y \sim Q}\left|f^c \circ h_u(Y) - g^c \circ h_v(Y)\right|^2\right)^{1/2} \\
		& = \rho((f,u),(g,v)).
	\end{align*}
	Moreover, for any $u \in \SS^{d-1}$, $\Phi(\cdot,u)$ achieves the maximum over $\sC$ at a unique $f_u \in \sC$, since under the assumption of \ref{CC}, the Kantorovich potential corresponding to $P_u$ and $Q_u$ is unique.

	We now apply the functional delta method to \cref{prop:weak} with the supremum function $\iota.$
	By Kantorovich duality, $\iota(\Phi_{(P, Q)})(u) = W_p^p(P_u, Q_u)$ and $\iota(\Phi_{(P_n, Q_n)})(u) = W_p^p(P_{nu}, Q_{nu})$
	\cref{prop:hadamard} implies that $\iota$ is Hadamard directionally differentiable and the derivative of $\iota$ at $\Phi_{(P,Q)}$ in the direction $\BB$ is given by \begin{equation*}
    		\iota'_\Phi(\BB_{(P,Q)})(u) = \BB_{(P,Q)}(f_u, u).
	\end{equation*}	
	Hence,
	\begin{equation} \label{eq:19}
    		\sqrt{n}\left(\iota(\Phi_{(P_n,Q_n)}) - \iota(\Phi_{(P,Q)})\right) \rightsquigarrow   \iota'_{\Phi_{(P,Q)}}(\BB)(\cdot) = \BB(f_\cdot, \cdot)
    		\quad \text{ in } \ell^{\infty}(\SS^{d-1}).
	\end{equation}
	This is a centered tight Gaussian process on $\SS^{d-1}$, and a direct computation shows that its covariance agrees with that of $\sqrt 2 \GG$, as desired.
	The case where $n \neq m$ follows similarly; the details are omitted.
\end{proof}

\subsubsection{Additional Lemmas}
The following lemma is included in the proof for Corollary 2.7.2 of~\cite{vw-ep}. For the sake of clarity, we state it here separately.
\begin{lemma} \label{lemma:2}
    Suppose that $P$ is a probability distribution on $\RR^d$, then the bracketing number of any class of functions $\sH$ with respect to $L_2(P)$ can be bounded above by the covering number with respect to the uniform norm. Explicitly, for $\varepsilon > 0$, we have \begin{equation} \label{eq:10}
    N_{[]}(2\varepsilon,\sH,L_2(P)) \leq N(\varepsilon,\sH,\|\cdot\|_{\infty}).
\end{equation}
\end{lemma} 
\begin{proof} The proof is inspired by that of Lemma 6.2 of~\cite{wellner-emp}. Take any $g \in \sH$ and suppose that it lies in some ball $B_{\|\cdot\|_{\infty}}(f,\varepsilon)$. Let $l = f - \varepsilon$, $u = f + \varepsilon$, then \begin{equation*}
    l - g = f - g - \varepsilon \leq \varepsilon - \varepsilon = 0,\, u - g = f - g + \varepsilon \geq -\varepsilon + \varepsilon = 0,
\end{equation*} and $\|l-u\|_{L^2(P)} = 2\varepsilon$. \end{proof}

\noindent The next lemma is included in the proof for Lemma 4.2 of~\cite{dm-cover}. Again we state and prove it here separately for completeness.
\begin{lemma} \label{lemma:3}
    Consider some set of composite functions $\sA_{\sH} := \{f \circ h:\,f \in \sA,\,h \in \sH\}$ where the functions in $\sA$ are $K$-Lipschitz. Then \begin{equation} \label{eq:11}
    N((K+1)\varepsilon,\sC_{\sF},\|\cdot\|_{\infty}) \leq N(\varepsilon,\sC_R,\|\cdot\|_{\infty}) \times N(\varepsilon,\sF,\|\cdot\|_{\infty}).
\end{equation}
\end{lemma} 
\begin{proof} 
Indeed, for any $g := f \circ h \in \sA_{\sH}$, there exists some $f_0 \in \sA$ and $h_0 \in \sH$ such that $\|f - f_0\|_{\infty},\,\|h - h_0\|_{\infty} < \varepsilon$, then \begin{align*}
    & \|g - f_0 \circ h_0\|_{\infty} \leq \|f \circ h - f \circ h_0\|_{\infty} + \|f \circ h_0 - f_0 \circ h_0\|_{\infty} \\
    & \leq K\|h - h_0\|_{\infty} + \|f - f_0\|_{\infty} < (K+1)\varepsilon.
\end{align*}
\end{proof}

\begin{lemma} \label{lemma:4}
	The covering entropy of $\sF^C$  with respect to $L^{\infty}$ is upper bounded by that of $\sF$. 
\end{lemma}
\begin{proof}
	Fix any $f^c \in \sF^C$. There exists some $f_0 \in \sF$ such that $\|f - f_0\|_{\infty} < \varepsilon$. Then \begin{align*}
		\|f^c - f_0^c\|_{\infty} & = \|\inf_x \left(|x-y|^p - f(x)\right) - \inf_x \left(|x-y|^p - f_0(x)\right)\|_{\infty} \\
		& =  \|\sup_x \left(f_0(x) - |x-y|^p\right) - \sup_x \left(|x-y|^p - f(x)\right)\|_{\infty} \\
		& \leq  \|\sup_x \left(|x-y|^p - f(x)) -(|x-y|^p - f_0(x)\right)\|_{\infty} \\
		& = \|\sup_x \left(f_0(x) - f(x)\right)\|_{\infty} \leq \|f_0 - f\|_{\infty} < \varepsilon\,.
	\end{align*}
	The conclusion follows immediately.
\end{proof}

\subsection{Proof of Theorem \ref{thm:2}}
\begin{proposition} \label{prop:brenier}
    Consider two probability distributions $P$ and $Q$ that satisfy \ref{CC} and have compact supports contained in the ball $\overline{B(0,R)}$ for some $R > 0$.
    Let $f_u, f_v \in \sC$ be the unique Kantorovich potentials for $(P_u, Q_u)$ and $(P_v, Q_v)$, respectively.
    Then there exists a constant $C_{R, p}$ depending on $R$ and $p$ such that
    \begin{equation}
    	\|f_u - f_v\|_\infty \leq C_{R, p} \|u - v\|_2^{p-1}.
    \end{equation}
    Moreover, if $P$ and $Q$ are discrete probability distributions on $\{x_1,\dots,x_N\},\{y_1,\dots,y_N\} \subset \overline{B(0,R)}$ respectively such that $P(x_i) = Q(y_i) = 1/N$ for $i = 1,\dots,N$, then the inequality above also holds.
\end{proposition}
\begin{proof} 
	By the representation of one-dimensional Wasserstein costs~\cite[see][Theorem 2.10]{bl-1d}, we have for any $u \in \SS^{d-1}$,
	\begin{equation*}
		W_p^p(P_u, Q_u) = \int_0^1 |P_u^{-1}(t) - Q_u^{-1}(t)|^p \dd t\,,
	\end{equation*}
	where, by abuse of notation, $P_u^{-1}$ and $Q_u^{-1}$ denote the inverses of the cumulative distribution functions of $P_u$ and $Q_u$, i.e.,
	\begin{equation*}
	    P_u^{-1}(t) = \inf \{x \in \RR: \PP_P\{X^\top u \leq x\} \geq t\}\,,
	\end{equation*}
	and analogously for $Q_u^{-1}$.
	Under \ref{CC}, note that these inverses satisfy 
	\begin{align*}
	    P_u^{-1}(t) \leq x & \iff t \leq P_u(x) \\
	    P_u^{-1}(t) \geq x & \iff t \geq P_u(x)\,,
	\end{align*}
	and likewise for $Q_u$ (this follows from the considerations in \cite[Section 2.1]{San15} combined with the fact that $P_u^{-1}$ is a right inverse for $P_u$ since the support of $P_u$ is connected).
	
	It follows from \cite[Theorem 1.17]{San15} that the derivative of any optimal Kantorovich potential $f_u$ must satisfy
	\begin{equation*}
		f_u'(x) = p |x - Q_u^{-1} \circ P_u(x)|^{p-2} (x - Q_u^{-1} \circ P_u(x))\,.
	\end{equation*}
	Note that $|x - Q_u^{-1} \circ P_u(x)| \leq 2R$, so that this expression is bounded by $p (2R)^{p-1}$.
	Therefore, if we define
	\begin{equation*}
		f_u(x) = \int_0^x p |x' - Q_u^{-1} \circ P_u(x')|^{p-2} (x' - Q_u^{-1} \circ P_u(x')) \dd x'\,,
	\end{equation*}
	then $f_u$ is $p (2R)^{p-1}$ Lipschitz, satisfies $f_u(0) = 0$, and is a Kantorovich potential, and under \ref{CC}, it must therefore be the unique optimal potential in $\sC$.
	
	If we define $g_u(x') = |x' - Q_u^{-1} \circ P_u(x')|^{p-2} (x' - Q_u^{-1} \circ P_u(x'))$, it follows that
	\begin{align*}
		\|f_u - f_v\|_\infty  & = \max_{x \in [-R, R]} \left|\int_0^x p (g_u(x') - g_v(x')) \dd x'\right| \\
		& \leq p \int_{-R}^R |g_u(x') - g_v(x')| \dd x'\,.
	\end{align*}
	The function $v \mapsto |v|^{p-2} v$ is $p - 1$-H\"older, with norm depending on $R$ and $p$~\cite[see, e.g][proof of Corollary 3]{manole2021sharp}.
	Letting $C_{R, p}$ denote a constant depending on $R$ and $p$ whose value may vary from line to line, we obtain
	\begin{align*}
		\|f_u - f_v\|_\infty & \leq C_{R, p} \int_{-R}^R |Q_u^{-1} \circ P_u(x') - Q_v^{-1} \circ P_v(x')|^{p-1} \dd x' \\ & C_{R, p} \int_{-R}^R |Q_u^{-1} \circ P_u(x') - Q_v^{-1} \circ P_u(x')|^{p-1} + |Q_v^{-1} \circ P_u(x') - Q_v^{-1} \circ P_v(x')|^{p-1} \dd x'
		\\ & \leq C_{R, p} (\|Q_u^{-1} \circ P_u - Q_v^{-1} \circ P_u\|_\infty^{p-1} + \int_{-R}^R |Q_v^{-1} \circ P_u(x') - Q_v^{-1} \circ P_v(x')|^{p-1} \dd x')
	\end{align*}
	For the first term, it suffices to note that $\|Q_u^{-1} - Q_v^{-1}\|_\infty$ is bounded.
	Indeed, \cite[equation (2.3)]{bl-1d} implies
	\begin{align*}
		\|Q_u^{-1} - Q_v^{-1}\|_\infty & = W_\infty(Q_u, Q_v) \\
		& \leq \|Y^\top u - Y^\top v\|_{L^\infty(Q)} \\
		& \leq R \|u - v\|_2\,,
	\end{align*}
	where the first inequality follows from the fact that $(Y^\top u, Y^\top v)$ with $Y \sim Q$ is a valid coupling of $Q_u$ and $Q_v$.
	Therefore $\|Q_u^{-1} \circ P_u - Q_v^{-1} \circ P_u\|_\infty \leq R \|u-v\|_2$.
	
	For the second term, we first derive an upper bound for the case $p = 2$. The idea is borrowed from the proof of \cite[Proposition 2.17]{San15}. Through computations, we have \begin{align*}
		& \int_{-R}^{R}|Q_v^{-1} \circ P_u(x') - Q_v^{-1} \circ P_v(x')|\dd x' \\
		& = \mathcal{L}^2\left(\left\{(x',y) \in [-R,R] \times [-R,R]: Q_v^{-1} \circ P_u(x') \leq y < Q_v^{-1}\circ P_v(x') \right.\right. \\
		& \left.\left. \hskip 15.5em \text{ or } Q_v^{-1} \circ P_v(x') \leq y < Q_v^{-1}\circ P_u(x')\right\}\right) \\
		& = \mathcal{L}^2\left(\left\{(x',y) \in [-R,R] \times [-R,R]: Q_v^{-1} \circ P_u(x') \leq y < Q_v^{-1}\circ P_v(x')\right\}\right) \\
		& \quad + \mathcal{L}^2\left(\left\{(x',y) \in [-R,R] \times [-R,R]: Q_v^{-1} \circ P_v(x') \leq y < Q_v^{-1}\circ P_u(x')\right\}\right)\,.
	\end{align*} 
	By Fubini's theorem along with the monotonicity of cumulative distribution functions, we have \begin{align*}
		& \mathcal{L}^2\left(\left\{(x',y) \in [-R,R] \times [-R,R]: Q_v^{-1} \circ P_u(x') \leq y < Q_v^{-1}\circ P_v(x')\right\}\right) \\
		& = \int_{-R}^{R}\mathcal{L}^1\left(\left\{x' \in [-R,R]: P_v^{-1}\circ Q_v(y) < x' \leq P_u^{-1}\circ Q_v(y) \right\}\right) \dd y\,.
	\end{align*}
	Similarly, \begin{align*}
		& \mathcal{L}^2\left(\left\{(x',y) \in [-R,R] \times [-R,R]: Q_v^{-1} \circ P_v(x') \leq y < Q_v^{-1}\circ P_u(x')\right\}\right) \\
		& = \int_{-R}^{R}\mathcal{L}^1\left(\left\{x' \in [-R,R]: P_u^{-1}\circ Q_v(y) < x' \leq P_v^{-1}\circ Q_v(y) \right\}\right) \dd y\,.
	\end{align*}
	Summing up the integrals, we obtain \begin{align*}
		 \int_{-R}^{R}|Q_v^{-1} \circ P_u(x') - Q_v^{-1} \circ P_v(x')|\dd x' & = \int_{-R}^{R} \left|P_u^{-1}\circ Q_v(y) - P_v^{-1}\circ Q_v(y)\right| \dd y \\
		 & \leq 2R\|P_u^{-1} - P_v^{-1}\|_{\infty}  \leq 2R^2\|u-v\|_2\,.
	\end{align*}
	When $p > 2$, we have \begin{align*}
		\int_{-R}^R |Q_v^{-1} \circ P_u(x') - Q_v^{-1} \circ P_v(x')|^{p-1} \dd x'&\leq (2R)^{p-2}\int_{-R}^R |Q_v^{-1} \circ P_u(x') - Q_v^{-1} \circ P_v(x')| \dd x' \\&\leq (2R)^p\|u-v\|_2\,.
	\end{align*}
	
	Finally, when $1 < p < 2$, again we consider $x'$ as a random variable of the uniform distribution $U'$ on $[-R,R]$, and then by Jensen's inequality \begin{align*}
		\int_{-R}^R |Q_v^{-1} \circ P_u(x') - Q_v^{-1} \circ P_v(x')|^{p-1} \dd x' & = 2R \frac{1}{2R}\int_{-R}^R |Q_v^{-1} \circ P_u(x') - Q_v^{-1} \circ P_v(x')|^{p-1} \dd x' \\
		& \leq 2R (\frac{1}{2R}\int_{-R}^R |Q_v^{-1} \circ P_u(x') - Q_v^{-1} \circ P_v(x')| \dd x')^{p-1} \\
		& \leq 2R (R\|u-v\|_2)^{p-1} \leq 2R^{p}\|u-v\|_2^{p-1}\,.
	\end{align*}
	Hence, combining the upper bounds for the first and second term, we obtain for $p > 1$, \begin{equation*}
		\|f_u - f_v\|_\infty \leq C_{R,p} (\|u-v\|_2^{p-1} + \|u-v\|_2^{(p-1) \vee 1}) \leq C_{R,p}\|u-v\|_2^{p-1}.
 	\end{equation*}
 	Now we turn to consider the discrete distributions. In this case, the inverse $P_u^{-1}$ only satisfies \begin{align*}
	    P_u^{-1}(t) \leq x & \iff t \leq P_u(x) \\
	    P_u^{-1}(t) > x & \iff t > P_u(x)\,,
	\end{align*}
	and likewise for $Q_u^{-1}$. 
	
	Fix any $y \in [-R,R]$, we have \begin{align*}
		& \mathcal{L}^1\left(\left\{x' \in [-R,R]: P_u(x') \leq Q_v(y) \text{ and } P_v(x') > Q_v(y) \right\}\right)  \\
		& = \mathcal{L}^1\left(\left\{x' \in [-R,R]: P_u(x') < Q_v(y) \text{ and } P_v(x') > Q_v(y)\right\}\right)  \\
		& \quad + \mathcal{L}^1\left(\left\{x' \in [-R,R]: P_u(x') = Q_v(y) \text{ and } P_v(x') > Q_v(y)\right\}\right) \\
		& \leq \mathcal{L}^1\left(\left\{x' \in [-R,R]: P_v^{-1} \circ Q_v(y) < x' < P_u^{-1} \circ Q_v(y)\right\}\right) \\
		& \quad + \mathcal{L}^1\left(\left\{x' \in [-R,R]: P_u^{-1}\circ Q_v(y) \leq x < P_u^{-1}\circ (Q_v(y) +1/N)\right\}\right) \\
		& = \mathcal{L}^1\left(\left\{x' \in [-R,R]: P_v^{-1} \circ Q_v(y) < x' < P_u^{-1} \circ (Q_v(y) +1/N) \right\}\right)\,.
	\end{align*}
	Analogously, \begin{align*}
		& \mathcal{L}^1\left(\left\{x' \in [-R,R]: P_v(x') \leq Q_v(y) \text{ and } P_u(x') > Q_v(y) \right\}\right) \\
		& \leq \mathcal{L}^1\left(\left\{x' \in [-R,R]: P_u^{-1} \circ Q_v(y) < x' < P_v^{-1} \circ (Q_v(y) +1/N) \right\}\right)\,.
	\end{align*}
	In addition, the two sets have no intersection. Therefore, \begin{align*}
		& \int_{-R}^{R}\left|Q_v^{-1} \circ P_u(x') - Q_v^{-1} \circ P_v(x')\right|\dd x'  \\
		& \leq \int_{-R}^{R} \left|P_u^{-1}\circ (Q_v(y)+1/N) - P_v^{-1}\circ Q_v(y) + P_u^{-1} \circ Q_v(y) -  P_v^{-1} \circ (Q_v(y) +1/N) \right| \dd y \\
		& \leq \int_{-R}^{R} \left|P_u^{-1}\circ (Q_v(y)+1/N) -  P_v^{-1} \circ (Q_v(y) +1/N) \right| + \left| P_u^{-1} \circ Q_v(y) - P_v^{-1}\circ Q_v(y)\right|  \dd y \\
		& \leq 4R\|P_u^{-1} - P_v^{-1}\|_{\infty} \leq 4R^2\|u-v\|_2\,.
	\end{align*}

\end{proof} 


\begin{proof}[Proof of \cref{thm:2}]
	We may assume without loss of generality that $m = n$ by discarding addition samples from either $P$ or $Q$, if necessary.
We define an estimator for $u,v \in \SS^{d-1}$ by setting\begin{equation} \label{eq:estimator}
    		\begin{split} \hat{\Sigma}_{u,v} = & \int f_{nu}(u^{\top}x)f_{nv}(v^{\top}x) \dd P_n(x) \\
    		& - \left( \int f_{nu}(u^{\top}x) \dd P_n(x)\right)\left(\int f_{nv}(v^{\top}x) \dd P_n(x)\right) \\
    		& + \int f_{nu}^c(u^{\top}y)f_{nv}^c(v^{\top}y) \dd Q_n(y) \\
    		& - \left(\int f_{nu}^c(u^{\top}y) \dd Q_n(y)\right)\left(\int f_{nv}^c(v^{\top}y) \dd Q_n(y)\right)\,,
    		\end{split}
		\end{equation}
	where $f_{n\cdot} \in \sC$ denotes a Kantorovich potential for $P_{n\cdot}$ and $Q_{n\cdot}$.
	We are done if we can show the convergence of the first term to the corresponding one in (\ref{eq:6}). The proof for the other three terms follows similar routine. 
	
	We first split the absolute difference between the objective quantities into two terms: \begin{align*}
		& \sup_{u,v \in \SS^{d-1}}\left|\int f_{nu}(u^{\top}x)f_{nv}(v^{\top}x) \dd P_n(x) - \int f_u(u^{\top}x)f_v(v^{\top}x) \dd P(x)\right| \\
		& \leq \frac{1}{n}\sup_{u,v \in \SS^{d-1}}\left|\sum_{i=1}^n f_{nu}(u^{\top}X_i)f_{nv}(v^{\top}X_i) - f_u(u^{\top}X_i)f_v(v^{\top}X_i)\right| \\ 
		& \quad + \sup_{u,v \in \SS^{d-1}}\left|\frac{1}{n}\sum_{i=1}^n \left( f_{u}(u^{\top}X_i)f_{v}(v^{\top}X_i) - \EE f_u(u^{\top}X_i)f_v(v^{\top}X_i)\right)\right|\,.
	\end{align*}
	For the second term, we notice that $P_n - P \rightsquigarrow 0$ in $\ell^{\infty}(\sC \times \SS^{d-1})$ and $H:\ell^{\infty}(\sC \times \SS^{d-1}) \to \RR$ defined by $H(g) = \sup_{u,v}|g(f_u,u)g(f_v,v)|$ is continuous. Therefore, the second term converges to zero in distribution and consequently in probability. Since all $f_{nu},f_u,\,u \in \SS^{d-1}$ are uniformly bounded by some constant depending only on $R$ and $p$, dominated convergence theorem yields convergence in mean.
	
	For the first term, we have
	\begin{align*}
    		& \EE \sup_{u,v \in \SS^{d-1}}\left|P_n\left((f_{nu}\circ h_u) (f_{nv}\circ h_v) - (f_u \circ h_u) (f_v \circ h_v)\right)\right|  \\
    		& \leq \EE \sup_{u,v \in \SS^{d-1}}\left|f_{nu}(u^{\top}X)f_{nv}(v^{\top}X) - f_u(u^{\top}X)f_v(v^{\top}X)\right| \\
		& \leq C_{R,p}\EE \sup_{u \in \SS^{d-1}} \left|f_{nu}(u^{\top}X) - f_u(u^{\top}X)\right|\,.
	\end{align*}	
	\cref{prop:brenier} implies that \begin{align} \label{eq:empuniform}
    		& \left\|\left|f_{nu}(u^{\top}x) - f_u(u^{\top}x)\right| - \left|f_{nv}(v^{\top}x) - f_v(v^{\top}x)\right|\right\|_{\infty} \nonumber \\
		& \leq \|(f_{nu}(u^{\top}x) - f_u(u^{\top}x)) - (f_{nv}(v^{\top}x) - f_v(v^{\top}x))\|_{\infty} \nonumber \\
    		& \leq C_{R,p}\|u-v\|_2^{p-1}.
	\end{align}
	Therefore, for any $\varepsilon > 0$, there exists some partition of $\SS^{d-1}$ such that for $n \in \NN^{*}$ large enough, \begin{align*}
    		& \SS^{d-1} \subseteq \mathop{\bigcup}\limits_{i=1}^{N_{\varepsilon}}B(u_i,\delta_{\varepsilon}),\\
    		& \left\|\left|f_{nu}(u^{\top}x) - f_u(u^{\top}x)\right| - \left|f_{nu_i}({u_i}^{\top}x) - f_{u_i}({u_i}^{\top}x)\right|\right\|_{\infty} < \frac{\varepsilon}{2},\,i = 1,\dots,N_{\varepsilon}, \\
    		& \EE \mathop{\sup}\limits_{i=1,\dots,N_{\varepsilon}}\left|f_{nu}(u^{\top}X) - f_u(u^{\top}X)\right| < \frac{\varepsilon}{2}.
	\end{align*}
	The last inequality follows from the the $P$-a.s. convergence of Kantorovich potentials by~\cite[Theorem 2.8]{bl-clt} with dominated convergence theorem applied to it and and the finiteness of terms taken over supremum.  
	
	Altogether, the first term may be bounded by arbitrarily small numbers when $n \to \infty$: \begin{equation*}
    		\EE\sup_{u \in \SS^{d-1}}\left|f_{nu}(u^{\top}X) - f_u(u^{\top}X)\right| \leq \EE \mathop{\sup}\limits_{i=1,\dots,N_{\varepsilon}}\left|f_{nu}(u^{\top}X) - f_u(u^{\top}X)\right| + \frac{\varepsilon}{2} < \varepsilon.
	\end{equation*}
	This completes the deduction of convergence in mean of the estimator \eqref{eq:estimator}. 
\end{proof}

\subsection{Proof of Theorem \ref{thm:3}}
\begin{proof}
	We define $H:\ell^{\infty}(\SS^{d-1}) \to \RR$ as $H(f) := \int_{\SS^{d-1}}f(\theta)\,d\sigma^d(\theta)$. $W_2^2(\mu_{\cdot},\nu_{\cdot})$ indeed belongs to $\ell^{\infty}(\SS^{d-1})$ since the Wasserstein distance between any one-dimensional projections of probability distributions $\mu$ and $\nu$ is bounded above by the one between $\mu$ and $\nu$ themselves. Besides, the integral over unit sphere with respect to uniform measure preserves the sup norm of the functions in $\ell^{\infty}(\SS^{d-1})$.
	
	By definition of weak convergence in $\ell^{\infty}(\SS^{d-1})$, the uniform CLT implies that \begin{equation*}
    \sqrt{n}(SW_p^p(P_n,Q_n) - SW_p^p(P,Q)) = \int_{\SS^{d-1}}\sqrt{n}(W_p^p(P_{n\theta},Q_{n\theta}) - W_p^p(P_{\theta},Q_{\theta}))\,d\sigma^d(\theta) \stackrel{\text{d}}{\to} S.
\end{equation*}
	For any $\omega$ in the probability space $\Omega$, $\int_{\SS^{d-1}}|\GG_{\theta}(\omega)|\,d\sigma^d(\theta) < \infty$. This can be easily deduced from the fact that $\GG$ has continuous sample paths a.s. In addition, $\GG:\SS^{d-1} \times \Omega \to \RR$ is jointly measurable and thus $\omega \mapsto S(\omega)$ is a random variable. Finally, $S$ is Gaussian due to the Riemann integrability of $\theta \mapsto \GG_{\theta}(\omega)$.
	
	Finally we compute the mean and variance of $S$. Trivially, $\EE S = 0$. In terms of the variance, we have \begin{align} \label{eq:variance}
    & \text{Var}(S) = \EE{\left(\int_{\SS^{d-1}}\GG(\theta) \dd \sigma(\theta)\right)}^2  = \int_{\SS^{d-1}}\int_{\SS^{d-1}}\EE(\GG(u)\GG(v)) \dd \sigma(u)\dd\sigma(v) \nonumber\\
    & = \int \left(\int_{\SS^{d-1}}f_{\theta}(\theta^{\top}x) \dd \sigma(\theta)\right)^2 \dd P(x) - \left(\int (\int_{\SS^{d-1}}f_{\theta}(\theta^{\top}x) \dd \sigma(\theta)) \dd P(x)\right)^2 \nonumber\\
    & + \int \left(\int_{\SS^{d-1}}f_{\theta}^c(\theta^{\top}y) \dd \sigma(\theta)\right)^2 \dd Q(x) - \left(\int (\int_{\SS^{d-1}}f_{\theta}^c(\theta^{\top}y) \dd \sigma(\theta)) \dd Q(y)\right)^2 \nonumber\\
    & = \text{Var}_{X \sim P}\left(\int_{\SS^{d-1}}f_{\theta}(\theta^{\top}X) \dd \sigma(\theta)\right) + \text{Var}_{Y \sim Q}\left(\int_{\SS^{d-1}}f_{\theta}^c(\theta^{\top}Y) \dd \sigma(\theta)\right)\,.
    \end{align}
	
\end{proof}

\begin{remark} \label{rmk:2}
	The variance of $S$ is identical with that of $Z_{(P,Q)}$ derived in Theorem 3 and Lemma 8 of~\cite{mbw-trim} with $\delta = 0$. The variance of $Z_{(P,Q)}$ can be reduced to \begin{equation*}
         \int_{\SS^{d-1}}\int_{\SS^{d-1}}\left(\int_0^1|P_u^{-1}(t)-Q_u^{-1}(t)|^p|P_v^{-1}(t)-Q_v^{-1}(t)|^p \dd t - W_p^p(P_u,Q_u)W_p^p(P_v,Q_v)\right) \dd \sigma(u)\dd \sigma(v).
    \end{equation*}
    Let $\pi_u$ and $\pi_v$ denote the optimal transport plans between $P_u,Q_u$ and $P_v,Q_v$ respectively. Letting $(X_u,Y_u) \sim \pi_u$ and $(X_v,Y_v) \sim \pi_v$, it follows that the expression above is equal to\begin{align*}
        & \int_{\SS^{d-1}}\int_{\SS^{d-1}} \text{cov }\left(|X_u-Y_u|^p,|X_v-Y_v|^2\right) \dd \sigma (u)\dd\sigma(v) \\
        & = \text{Var}\left(\int_{\SS^{d-1}}|X_{\theta}-Y_{\theta}|^p \dd\sigma(\theta)\right) \\
        & = \text{Var}\left(\int_{\SS^{d-1}}f_{\theta}(\theta^{\top}X) + f_{\theta}^c(\theta^{\top}Y) \dd\sigma(\theta)\right) \\
        & = \text{Var}_{X \sim P}\left(\int_{\SS^{d-1}}f_{\theta}(\theta^{\top}X)\dd \sigma(\theta)\right) + \text{Var}_{Y \sim Q}\left(\int_{\SS^{d-1}}f_{\theta}^c(\theta^{\top}Y) \dd \sigma(\theta)\right).
    \end{align*}
\end{remark} 

\section{Additional Experiments}
\subsection{Sliced Wasserstein Distance.}
 Consider the example in section \ref{simusw}. Instead of $p = 2$, we investigate the asymptotic behavior of the case $p = 1$. We first give an explicit representation of the theoretical limit of the example given in section \ref{simusw}. Then the unique $1$-Lipschitz function that achieves the $1$-Wasserstein distance between $P_\theta$ and $Q_\theta$ is $\phi_{0}^{\theta}(x) = -\text{sign}(a_{\theta})x$. Hence, we have \begin{equation*}
    \sqrt{n}\left(W_1(P_{n\cdot},Q_{n\cdot} - W_1(P_{\cdot},Q_{\cdot})\right) \rightsquigarrow \GG,
\end{equation*} where $\GG$ is the mean-zero Gaussian process indexed by $\SS^2$ with covariance functions \begin{equation*}
    \EE\GG(u)\GG(v) = \frac{2}{3}\text{sign}(a_u)\text{sign}(a_v)\langle u,v \rangle.
\end{equation*}
It follows from Theorem \ref{thm:3} that the limiting distribution of the empirical $1$-Wasserstein distance is the centered Gaussian $S$ with variance \begin{equation*}
    \text{Var}(S) = \frac{2}{3} \int_{\SS^2}\int_{\SS^2} \text{sign}(a_u)\text{sign}(a_v)\langle u,v \rangle\,d\sigma^3(u)d\sigma^3(v) \approx 0.164.
\end{equation*}
We sample i.i.d. observations $X_1,\dots,X_n \sim P$ and $Y_1,\dots,Y_n \sim Q$ with size $n = 50,100,500$. This process is repeated $500$ times. We then compare the finite distributions of $1$-Wasserstein distance with the theoretical limit given in section \ref{appsw}. We demonstrate the results using kernel density estimators in Figure \ref{fig:1}. We see that the finite-sample empirical distribution gets closer to the limiting Gaussian distribution in \ref{eq:29} as the sample size $n$ increases. \\
 In addition, we simulate the re-scaled plug-in bootstrap approximations by sampling $n = 1000$ observations of $P$ and $Q$. Fix some empirical SW $\sqrt{n}SW_2^2(P_n,Q_n)$, we generate $B = 500$ replications of $\sqrt{l}(SW_1(\hat{P}_n^*,\hat{Q}_n^*) - SW_1(P_n,Q_n))$. The distributions of the replications with various replacement numbers $l$, compared with the finite-sample empirical distribution and the theoretical limit, are shown in Figure \ref{fig:2}. We observe that the naive bootstrap ($l = n$) better approximates the finite sample distribution compared to fewer replacements ($l = n^{1/2},n^{3/4}$). This is consistent with the observation of inference on finite spaces.~\cite{oi-infproj} \\
\begin{figure} [!ht]
    \centering
    \begin{subfigure}{0.3\textwidth}
        \centering
        \includegraphics[width=1.\textwidth]{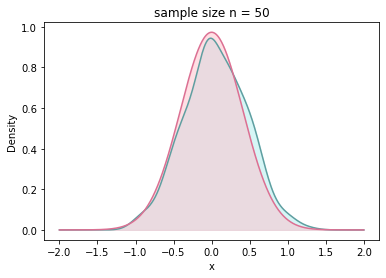} 
    \end{subfigure}
    \begin{subfigure}{0.3\textwidth}
        \centering
        \includegraphics[width=1.\textwidth]{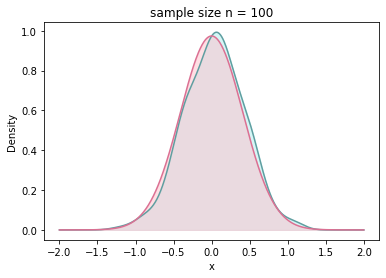} 
    \end{subfigure}
    \begin{subfigure}{0.3\textwidth}
        \centering
        \includegraphics[width=1.\textwidth]{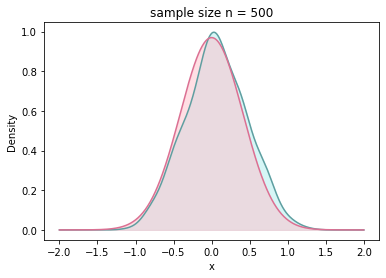} 
    \end{subfigure}
    \caption{Comparison of the finite sample density (pale turquoise) and the limit distribution of the empirical sliced distance (pink).}
    \label{fig:1}
\end{figure}
\begin{figure}[!ht]
    \centering
    \begin{subfigure}{0.33\textwidth}
        \centering
        \includegraphics[width=1.\textwidth]{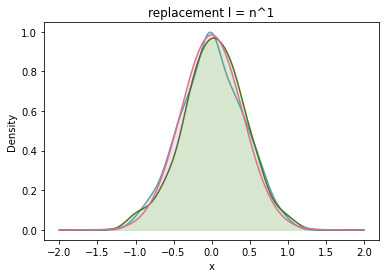} 
    \end{subfigure}
    \begin{subfigure}{0.33\textwidth}
        \centering
        \includegraphics[width=1.\textwidth]{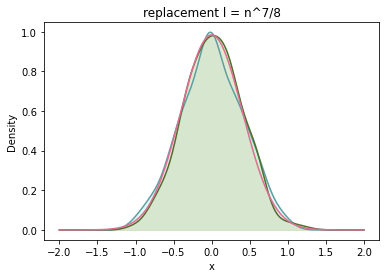} 
    \end{subfigure}
    \begin{subfigure}{0.33\textwidth}
        \centering
        \includegraphics[width=1.\textwidth]{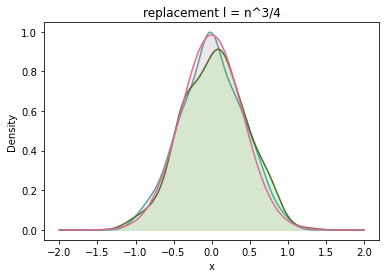} 
    \end{subfigure}
    \begin{subfigure}{0.33\textwidth}
        \centering
        \includegraphics[width=1.\textwidth]{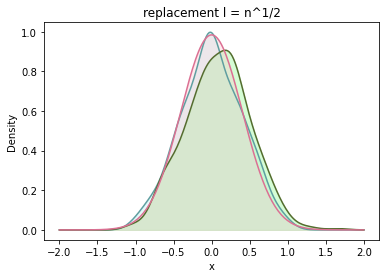} 
    \end{subfigure}
    \caption{Bootstrap for the empirical sliced distance. Illustration of the re-scaled plug-in bootstrap approximation ($n = 1000$) with replacement $l \in \{n, n^{7/8}, n^{4/5}, n^{1/2}\}$. Finite bootstrap densities (pale green) are compared to the corresponding finite sample density (pale turquoise) and the limit distribution (pink).}
    \label{fig:2}
\end{figure}


\subsection{Max-Sliced Wasserstein Distance}

 Consider the example in section \ref{simuwpp}. Again we estimate the distributional limit of the empirical distributions of $1$-Wasserstein distance but with $a = 2$. The unique $1$-Lipschitz function that achieves $1$-Wasserstein distance between $P_{e_1}$ and $Q_{e_1}$ or equivalently $P_{-e_1}$ and $Q_{-e_1}$ is $\phi_0^{e_1}(x) = -|x|$. Consequently, the theoretical limit in this case is the mean-zero Gaussian with variance \begin{equation*}
    \text{Var}(\GG_{\pm e_1}) = \frac{1}{2}\int_{-1}^1 x^2\,dx - \left(\frac{1}{2}\int_{-1}^1 -|x|\,dx\right)^2 + \frac{1}{4}\int_{-2}^2 y^2\,dy - \left(\frac{1}{4}\int_{-2}^2 |y|\,dy\right)^2 = \frac{5}{12}.
\end{equation*} 
The plots of comparison between the theoretical limit and the finite sample distributions of $n = 100,500,1000$ each of which is repeated $1000$ times are given in the bottom part of Figure \ref{fig:3}. The simulation of bootstrap is plotted in the last $4$ pictures of Figure \ref{fig:4}. The naive bootstrap ($l = n$) better approximates the finite sample distribution compared to fewer replacements ($l = n^{1/2},n^{3/4}$). \\

\begin{figure}[!ht]
    \centering
     \begin{subfigure}{0.3\textwidth}
        \centering
        \includegraphics[width=1.\textwidth]{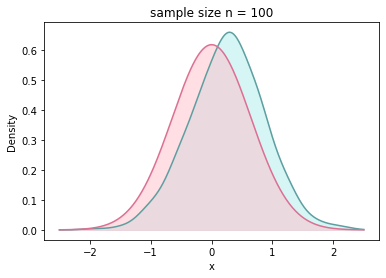}
    \end{subfigure}
    \begin{subfigure}{0.3\textwidth}
        \centering
        \includegraphics[width=1.\textwidth]{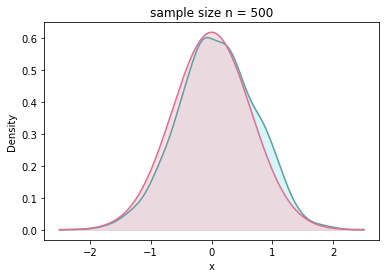} 
    \end{subfigure}
    \begin{subfigure}{0.3\textwidth}
        \centering
        \includegraphics[width=1.\textwidth]{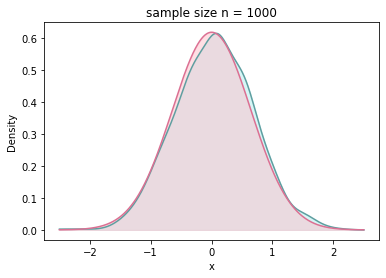} 
    \end{subfigure}
        \caption{Comparison of the finite sample density (pale turquoise) and the limit distribution of the empirical one-dimensional WPP (pink).}
    \label{fig:3}
\end{figure}
\begin{figure}[!ht]
    \centering
    \begin{subfigure}{0.33\textwidth}
        \centering
        \includegraphics[width=1.\textwidth]{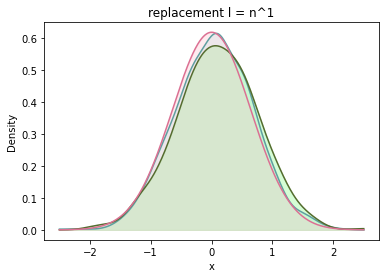} 
    \end{subfigure}
    \begin{subfigure}{0.33\textwidth}
        \centering
        \includegraphics[width=1.\textwidth]{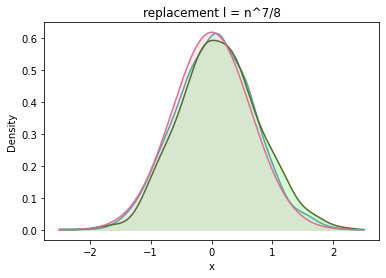}  
    \end{subfigure}
    \begin{subfigure}{0.33\textwidth}
        \centering
        \includegraphics[width=1.\textwidth]{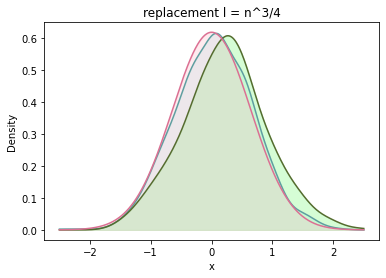}  
    \end{subfigure}
    \begin{subfigure}{0.33\textwidth}
        \centering
        \includegraphics[width=1.\textwidth]{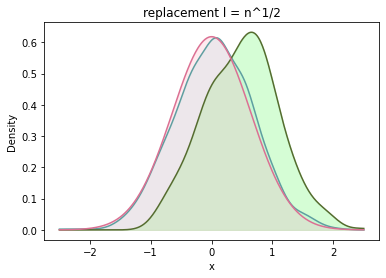}  
    \end{subfigure}
    \caption{Bootstrap for the empirical one-dimensional WPP. Illustration of the re-scaled plug-in bootstrap approximation ($n = 1000$) with replacement $l \in \{n, n^{7/8}, n^{4/5}, n^{1/2}\}$. Finite bootstrap densities (pale green) are compared to the corresponding finite sample density (pale turquoise) and the limit distribution (pink).}
    \label{fig:4}
\end{figure}


\subsection{Amplitude}

While it is hard to illustrate the uniform CLT, we may test it using certain functions. Instead of supremum, we consider the amplitude, i.e. $\text{amp}(f):= \sup f - \inf f$ for $f \in \ell^{\infty}(\SS^{d-1})$. The amplitude functional defined on $\SS^2$ is proven to be Hadamard directionally differentiable.~\cite{crc-diff} \\

\begin{figure}
    \centering
    \begin{subfigure}{0.33\textwidth}
        \centering
        \includegraphics[width=1.\textwidth]{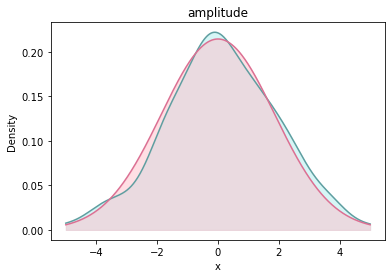}
    \end{subfigure}
    \begin{subfigure}{0.33\textwidth}
        \centering
        \includegraphics[width=1.\textwidth]{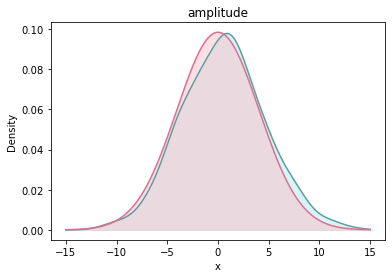} 
    \end{subfigure}
    \caption{Comparison of the finite sample density (pale turquoise) and the limit distribution of the empirical PRW (pink).
    \\Left: $\mathcal{P}\sim U(\{\frac{x^2}{4} + 4y^2 + z^2 = 1\})$, Right: $\mathcal{P}\sim U(\{\frac{x^2}{4} + \frac{y^2}{4} + \frac{z^2}{16} = 1\})$.\,\,\,$\mathcal{Q} \sim U(\SS^2)$.}
    \label{fig:10}
\end{figure}

\noindent Consider two distributions $P$ which is uniform on the surface of the ellipsoid $\{x^2/4 + 4y^2 + z^2 = 1\}$ and $\mathcal{Q}$ the uniform distribution on $\SS^2$. Apply delta method to $P,Q$ and $\text{amp}$, we get \begin{equation*}
    \sqrt{n}\left(\text{amp}(W_2^2(P_{n\cdot},Q_{n\cdot})) - 5/4\right) \stackrel{\text{d}}{\to} \GG((1,0,0)).
\end{equation*} We simulate the density of the amplitude of empirical Wasserstein distances of 1d projections. The finite sample density generated by $n = 600$ samples and theoretical limit are given in Figure \ref{fig:10}. \\

\noindent Let $P$ be uniform on $\{x^2/4 + y^2/4 + z^2/16 = 1\}$ and keep $Q$ unchanged. Then \begin{equation*}
    \sqrt{n}\left(\text{amp}(W_2^2(P_{n\cdot},Q_{n\cdot})) - 8/3\right) \stackrel{\text{d}}{\to} \GG((0,0,1)) - \GG((0,1,0)).
\end{equation*} We generate $n = 600$ samples according to $P$ and $Q$  and the result is also shown in Figure \ref{fig:10}.

\noindent Both finite sample densities indeed converge to the theoretical Gaussian limits as the sample size increases. This again serves as evidence of the uniform CLT.

\end{document}